\newcommand{\Label}[1]{\label{#1}\hspace{.3cm}\fbox{\rm
#1}\hspace{.3cm}}
\renewcommand{\Label}{\label} 
\newcommand{\marker}[1]{\mbox{\Huge$\bullet$}#1\mbox{\Huge$\bullet$}}
\renewcommand{\marker}[1]{#1}
\newcommand{\Comments}[1]{\\ \fbox{\fbox{\rm #1}}\\}
\renewcommand{\Comments}[1]{\mbox{}}
\documentclass[11pt]{article}
 \usepackage{amsmath,amsfonts,amssymb}
 \usepackage{amsfonts}
 \usepackage{amssymb}
\newcommand{\ls}[1]
   {\dimen0=\fontdimen6\the\font \lineskip=#1\dimen0
 \advance\lineskip.5\fontdimen5\the\font \advance\lineskip-\dimen0
 \lineskiplimit=.9\lineskip \baselineskip=\lineskip
 \advance\baselineskip\dimen0 \normallineskip\lineskip
 \normallineskiplimit\lineskiplimit \normalbaselineskip\baselineskip
 \ignorespaces }
 \numberwithin{equation}{section}
 \newtheorem{theorem}{Theorem}[section]
 
 \newtheorem{corollary}[theorem]{Corollary}
 \newtheorem{lemma}[theorem]{Lemma}
 
 \newtheorem{proposition}[theorem]{Proposition}
 \newtheorem{definition}[theorem]{Definition}
 \newtheorem{remark}[theorem]{Remark}
 \newtheorem{remarks}[theorem]{Remarks}

 \newcommand{\beq}{\begin{equation}}
 \newcommand{\eeq}{\end{equation}}
 \newcommand{\beqa}{\begin{eqnarray}}
 \newcommand{\eeqa}{\end{eqnarray}}
 \newcommand{\beqann}{\begin{eqnarray*}}
 \newcommand{\eeqann}{\end{eqnarray*}}

 \newcommand{\pf}{\noindent \mbox{{\bf Proof}: }}
 \def\squarebox#1{\hbox to #1{\hfill\vbox to #1{\vfill}}}
 \newcommand{\qed}{\hspace*{\fill}
    \vbox{\hrule\hbox{\vrule\squarebox{.667em}\vrule}\hrule}\smallskip}
 \newcommand{\req}[1]{(\ref{#1})}
 \newcommand{\lip}{\langle}
 \newcommand{\rip}{\rangle}

 \newcommand{\noi}{\noindent}

 \newcommand{\Law}{{\operatorname{Law\ }}}
 
 \newcommand{\grad}[1]{\nabla #1}
 \newcommand{\pd}[2]{\frac{\partial #1}{\partial #2}}
 \newcommand{\one}{{\frac{1}{n}}}%
 \newcommand{\half}{{\frac{1}{2}}}%

 \newcommand{\sfrac}[2]{\mbox{$\frac{#1}{#2}$}}

 \newcommand{\hilbert}{\bigcirc\kern -0.8em{\rm\scriptstyle {H}\;}}
 \newcommand{\la}{\lambda}
 \newcommand{\laa}{^{[\lambda]}}
 \newcommand{\eps}{\varepsilon}
 \newcommand{\om}{\omega}
 
 \newcommand{\vph}{\varphi}
 
 \newcommand{\bfcdot}{{\boldsymbol \cdot}}
 \newcommand{\Ds}{\Delta_s}
 \newcommand{\DD}{\mathbb{D}}
 \newcommand{\calE}{\mathcal E}
 \newcommand{\FF}{\calF}
 \newcommand{\FFF}{F}
 
 \newcommand{\calF}{{\mathcal F}}%
 \newcommand{\calL}{{\mathcal L}}%
 \newcommand{\calP}{{\mathcal P}}%
 \newcommand{\calS}{{\mathcal S}}%
 \newcommand{\upn}{^{(n)}}
 \newcommand{\wteta}{\widetilde{\eta}}
 \newcommand{\wtv}{\widetilde{v}}
 
 \newcommand{\wtT}{\widetilde{T}}

 \newcommand{\wtron}{\widetilde{\rho}^{(n)}}
 \newcommand{\wtTn}{\wtT^{(n)}}
 \newcommand{\wtvn}{\wtv^{(n)}}
 \newcommand{\Nat}{\mathbb{N}}
 \newcommand{\Reals}{\mathbb{R}}%
 \newcommand{\imb}{\,\mbox{\Large $\subset$}\!\!\!\!
               \raisebox{-.55ex}{\mbox{\footnotesize $\to$}}\,}
 \newcommand{\pair}[4]{\raisebox{-1.2ex}{\mbox{\tiny$#1$}}\!
       \langle#2,#3\rangle\!\raisebox{-1.2ex}{\mbox{\tiny $#4$}}}
 \newcommand{\clf}[2]{\pair{W}{#2}{#1}{W^*}}
 \newcommand{\clfs}[2]{\pair{W^*}{#2}{#1}{\Wss}}
 \newcommand{\clfss}[2]{\pair{W^*}{#1}{#2}{\Wss}}
 \newcommand{\clfsss}[2]{\pair{\Wss}{#1}{#2}{W^*}}
 
 \newcommand{\clfY}[2]{\pair{Y}{#1}{#2}{Y^*}}
 \newcommand{\conde}[2]
           {E\left(#1\,\raisebox{-.15cm}{\rule{.1mm}{5mm}}\,#2\right)}
 \newcommand{\EFn}[1]{\conde{#1}{\calF_n}}
 \newcommand{\dom}{{\rm dom}}
 \newcommand{\ddom}{{\rm \bf{dom}}}
 
 \newcommand{\domdel}[1]{\dom_{#1}\delta}
 \newcommand{\ddomdel}[1]{\ddom_{#1}\ddelta}
 \newcommand{\Domdel}[1]{{\ddomdel{#1}}}

 \newcommand{\ddelta}{\delta\hspace{-.2cm}\delta}
 \newcommand{\domo}[1]{\dom^{\mbox{\tiny $0$}}_{#1}\delta}
 
 \newcommand{\Dp}[1]{\mathbb{D}_{p,1}^#1}
 \newcommand{\Le}[1]{L^{#1}_{\rm e}(\mu;H)}
 \newcommand{\Lp}{L^p(\mu)}
 \newcommand{\LpW}{L^p(\mu;W)}
 \newcommand{\Wss}{W^{**}}
 \newcommand{\wn}[1]{\widehat{#1}^{^n}}
 \newcommand{\AAA}{\mathbf A}
 \newcommand{\BBB}{\mathbf B}
 \newcommand{\KKK}{\mathbf K}
 \newcommand{\TTT}{\mathbf T}
 
\begin{document}
 \ls{1.3}
   \title{The divergence of Banach space valued random variables on Wiener space}
   \author{E. Mayer-Wolf
           \footnote{Department of Mathematics, Technion I.I.T., Haifa, Israel}
       \ \ and M. Zakai\footnote{Department of Electrical Engineering,
                                       Technion I.I.T., Haifa, Israel}}
 \date{}
 \maketitle \vspace{.2cm}

  \hspace*{2.5cm}   
\fbox{\LARGE Corrigendum posted in arXiv:0710.4483}\vspace{.35cm}

 \begin{abstract}
 The domain of definition of the divergence operator $\delta$ on an abstract
 Wiener space $(W, H, \mu)$ is extended to include $W$--valued and
 $W\!\otimes\!W$-- valued ``integrands''. The main properties and
 characterizations of this extension are derived and it is shown that in some
 sense the added elements in\ $\delta$'s extended domain have divergence zero.
 These results are then applied to the analysis of quasiinvariant flows
 induced by $W$-valued vector fields and, among other results, it turns out
 that these divergence-free vector fields ``are responsible" for generating
 measure preserving flows.
 \end{abstract}
%
 \vfill \noindent{\sc Key Words:} Abstract Wiener Space, Divergence, Flows.

\newpage
 \section{Introduction}

 The classical Malliavin calculus is based on the notions of
 gradient and divergence operators in a Hilbert space setting.  The
 gradient operation deals traditionally with the directional
 derivative of real or Hilbert space valued random variables in the
 direction of elements of the Cameron-Martin space and the
 divergence operator is introduced by duality with respect to the
 gradient operator.  This setup, centered around separable Hilbert
 spaces, has proved to be a powerful tool in solving many problems.
 However, it needed to be extended in several cases,  whether in
 stochastic analysis on manifolds (\cite{D92}, \cite{D95},  \cite{CM96}),
 in the analysis on abstract Wiener spaces, (\cite{FP}, \cite{MD}, \cite{S94},
 \cite{P})  or in considerations associated with extending the Malliavin
 calculus  to include measure preserving transformations of the Wiener path
 (\cite{HUZ}, \cite{Z}).

 Let $\{W,H,\mu\}$ be an abstract Wiener space \marker{(cf. next section)}, let
 $\{e_i\}_{i\in\mathbb{N}}$ be a smooth ONB in $H$
 and $\{\eta_i\}_{i\in\mathbb{N}}$ a sequence of i.i.d.\ $N(0,1)$
 random variables on the Wiener space. Then, by the Ito-Nisio theorem
 (\cite{LT})\ $Y_n\!=\!\sum_1^n \eta_i e_i$
 converges in $W$ and the limit, say $Y$, is a measure preserving
 transformation on $\{W,\mu\}$.  There is no reason to expect that
 in general the difference $Y(\om)\!-\!\om$ will be $H$ valued (indeed,
 $Y(\omega)\!=\!-\om$ is one such counterexample).  Moreover, for a
 collection\ $Y_t = \sum_1^\infty \eta_i (t,\om) e_i$\ of such measure
 preserving transformations, $(d Y_t/dt)_{t=0}$ even if it exists,
 need not be $H$-valued. Consequently, the analysis of measure
 invariance (and related) flows on Wiener space requires the study
 of $W$-valued, rather than only Cameron-Martin, vector fields
 (\cite{CC}, \cite{HUZ}, \cite{CM00}, \cite{Z}).

 In this paper we (a) extend the domain of definition of the divergence
 operator to include Banach space valued random variables and derive the main
 properties and representation of this extension and (b) apply the results of
 the first part to the analysis of flows on Wiener space.

 In the next section we first summarize the background and notation for
 later reference. Differentiation of random variables is generalized by
 stipulating differentiability subspaces other than $H$, smaller or larger,
 yielding Sobolev spaces which respectively contain or are contained in
 the standard ones $\mathbb{D}_{p,1}$.\\
 \indent
 In Section~\ref{divsection} we extend the domain of definition of the
 divergence from $H$-valued to appropriate $W$-valued, and even\ $\Wss$--valued,
 random variables.  The main properties of the extended divergence are derived
 and Shigekawa's decomposition \cite{S} of the domain of this divergence into
 \textit{exact\/} and \textit{divergence-free\/} subspaces of ``integrands'' is
 shown to hold in this generalized setup as well. In fact it turns out that it
 is the class of divergence-free integrands which is extended but not the class
 of exact integrands.\\
 \indent
 The classical divergence also operates on $H^{\otimes 2}$, the Hilbert-Schmidt
 operators on $H$. Section~\ref{divsection} also contains its extension in this
 case to random operators from\ $W^*$\ to an arbitrary Banach space\ $Y$. This
 construction is then applied in Section~\ref{domo} for\ $Y\!=\!\Wss$\ (in this
 case operators from\ $W^*$\ to\ $\Wss$\ can be seen as bilinear forms on\ $W^*$)
 to derive the representation of any divergence-free integrand as the divergence
 of a random antisymmetric bilinear form on $W^*$.
 In \cite{S}, Shigekawa constructed a general setup for $H$-valued differential
 forms on Wiener space and derived the $H$ Hodge-Kodaira theory for this setup.
 Our results constitute an extension of this theory beyond\ $H$,\ restricted
 to forms of order~1 ($\Wss$--valued random variables) and of order~2 (random
 bilinear forms on $W^*$\marker{, not necessarily antisymmetric}).

 Section~\ref{flowsection} starts with an introduction to measure preserving
 transformation of Wiener space.  The results of Sections~\ref{divsection}
 and~\ref{domo} are applied in Section~\ref{flowsection} to derive new results
 concerning flows generated by $W$--valued vector fields, extending the results
 of~\cite{HUZ} and~\cite{Z} on measure preserving flows to general flows.

 Section~\ref{addrem} deals with (a)~The notion of adapted $W$-valued vector
 fields and conditions under which the flows they generate are adapted and
 (b)~The relation between the flow equation of Section~\ref{flowsection} and
 a class of scalar valued partial differential equations motivated by the
 non-random case introduced by P-L.~Lions in \cite{L}.

 In some of the results presented in this paper it is required that a
 $W$-valued random variable, say $u(\om)$, or a collection of such r.v.'s,
 have the ``representability" property that for some orthonormal basis
 $\{e_i\}_{i\in\mathbb{N}}$ of $H$ whose elements are in $W^*$
 \begin{equation} \label{newstar}
  \left\| u-\sum_{i=1}^n\:
  \clf{e_i}{u} e_i\right\|_W\: \underset{n\to\infty}{\longrightarrow}\; 0
  \hspace{1.5cm}{\rm a.s.\ \ \ or\ in\ \ }L^p(\mu)
 \end{equation}
 (the dual use of $e_i$ both as an element of $W$ and of $W^*$ will be
 further clarified later).\\
 The representation~\req{newstar} of $u$ will obviously hold if
 $\{e_i\}_{i\in\mathbb{N}}$ is a Schauder basis of $W$, but it might still
 be valid even if $W$ does not possess such a basis.
 Note that if $u$ is representable in this sense, and if $T:\,W\to W$ is
 measurable with $T^*\mu \ll \mu$, then so is $u\circ T$.
 Finding appropriate conditions under which a $W$-valued r.v $u$ is
 representable as in~\req{newstar} seems to be delicate.
 \bigskip

 \noindent
 {\bf Acknowledgement:}\   We wish to thank A.S. \"{U}st\"{u}nel for some
 uselful discussions and for calling our attention to reference~\cite{L}.
 \section{Preliminaries}
 In this section we first recall some notions of stochastic analysis in
 abstract Wiener space,   as well as the Ornstein--Uhlenbeck semigroup and
 its generator, the number operator. In the second part we generalize the
 notion of subspaces of differentiability to other than $H$.

 \subsection{Notation and Generalities} \Label{calL}
  The basic object in this paper will be an (infinite dimensional) abstract
  Wiener space\ $(W,H,\mu)$,\ \marker{based on a separable Banach space\ $W$\
  with a densely embedded Hilbert space\ $H$, and a Gaussian measure $\mu$ on\
  $W$ under which each $l\!\in\!W^*$ becomes an $N(0,|h|^2_H)$ random variable.
  The embeddings\ $i:H\to W$ and $i^*:W^*\to H$ will not always be written
  explicitly; thus, for example, an element $e\in W^*$ will also be considered
  to be an element in $H$ or in $W$, the distinction being clear from the
  context, as for example\ in \eqref{newstar}}.

  In $H$, the inner product is denoted by\ $(\cdot,\cdot)_{_H}$\ \marker{and
  the notation\ $|\cdot|_H$ for the norm in $H$ has already been used in the
  previous paragraph.}
  An orthonormal basis (ONB)\ \ $\calE\!=\!\{e_i\}_{i=1}^\infty$ of $H$ will
  be said to be {\bf smooth} if $e_i\!\in\!W^*$ for all $i$.\
  The norms in $W$ and $W^*$ are\ $\|\cdot\|_{_W}$\ and\ $\|\cdot\|_{_{W^*}}$\
  respectively, while the natural pairing between \ $l\!\in\!W^*$
  and $w\!\in\!W$\ \ (resp. between\ $w^{**}\!\in\!\Wss$\ and\ $l\!\in\!W^*$)
  \ \ is denoted $\clf{l}{w}$\ \ (resp.\ $\clfs{w^{**}}{l}$).\ \ Any of these
  subscripts may be omitted if no confusion arises.

  We recall the canonical zero--mean Gaussian field\
  $\{\delta h,\ h\!\in\!H\}$ whose correlation is given by $H$'s
  inner product. In particular, $\delta l(\om)\!=\!\clf{l}{\om}\ \ $a.s\ \
  for every\ $l\!\in\!W^*$.\ \
  For $1\!\le\!p\!\le\!\infty$,\ \ $L^p(\mu)$ or $L^p(W,\mu)$
 will denote $L^p(W,\calF,\mu)$
  where $\calF:=\sigma(\delta h,h\in H)$,\ the sigma--algebra generated by the
  canonical Gaussian field. The same applies to $L^p(\mu;Y)$ for any other
  Banach space $Y$.

  The space of bounded linear operators from a Banach space\ $X$\ to a Banach
  space\ $Y$\ is denoted\ $L(X,Y)$ equipped with the  operator norm\
  $\|A\|_{_{L(X,Y)}}\!=\!\sup\{\|Ax\|_{_Y},\ \|x\|_{_X}\!\le\!1\}$\ and
  $L(X)\!:=\!L(X,X).$\ The space of bilinear forms on a Banach
  space\ $X$\ is denoted\ $M_2(X)$\ and is equipped with the norm\
  $\|T\|_{_{M_2(X)}}\!=\!\sup\left\{|T(x,x')|,\ x,x'\in X,\
                                  \|x\|_{_X}\!=\!\|x'\|_{_X}\!=\!1\right\}$.

 The reader is assumed to be familiar with the basic notions of the
 Malliavin calculus, i.e., the gradient $\grad$ and the divergence $\delta$
 applied to the Sobolev spaces $\DD_{p,k}$\ \ \
 ($\nabla:\DD_{p,k}\to\DD_{p,k-1}(H)$ and
  $\delta:\DD_{p,k}(H)\to\DD_{p,k-1}$).\
%
%
%
 We will however be somewhat more explicit about the Ornstein-Uhlenbeck
 semigroup and conclude this subsection with a summary of some of its
 associated facts as needed in later sections (cf., e.g., \cite{S94},
\cite{U}, \cite{NZ}).

 Let $(\widetilde W,H,\mu)$ be an independent copy of $(W,H,\mu)$ and
 $f\!\in\!L^p(W,\mu)$. The Ornstein-Uhlenbeck semigroup on $L^p(\mu),\ p\ge 1$,
 is defined by the Mehler formula
 \begin{equation} \Label{starK}
  T_tf(\om)=E_{_W}f\left(e^{-t}\om+\sqrt{1-e^{-2t}}\:\widetilde\om\right)
 \end{equation}
 where $E_W$ denotes the conditional expectation conditioned on $W$. The
 family $\{T_t\}_{t\ge 0}$\ is a contraction, self-adjoint semigroup,
 whose infinitesimal generator\ $-\calL$\ \ satisfies
 \begin{equation} \label{starstarK}
  (1+ \calL)^{-\beta} = \frac{1}{\Gamma(\alpha)} \int_0^\infty t^{\beta-1}
                                                             e^{-t} T_t dt\,.
 \end{equation}
 Consequently $(1+\calL)^{-\beta}$, $\beta >0$, is a bounded operator on
 $L_p(\mu)$\ for every $\beta\!>\!0$. Moreover,
 \begin{itemize}
  \item[(i)] Since $T_t$ is self-adjoint, so are $(1\!+\!\calL)^{-\beta},\
   \beta>0$,\ and\ $\calL$;
  \item[(ii)] $\calL = \delta \circ \grad$ on $\DD_{p,2}$;
  \item[(iii)] $\grad(1+\calL)^{-\frac{1}{2}}$ is a bounded linear operator
   from $L^p(\mu)$ to $L^p(\mu;H)$ for any $p\in(1,\infty)$;
  \item[(iv)] $(1\!+\!\calL)^\beta\grad f=\grad\calL^\beta f$
   for all real\ $\beta$\ and every\ $f\!\in\!\DD_{p,1}$\ with\ $Ef\!=\!0$.
  \end{itemize}
  \bigskip

 The definition of $T_t$ can be extended to\ $f$'s taking values in a
 separable Banach space\ $Y$\ for which\ $E\|f\|_Y^p\!<\!\infty$\ (i.e. to
 $L^p(\mu;Y)$)\ in which case the expectation in~\req{starK} is defined as
 a Bochner integral. Formula~\req{starstarK} and the boundedness of\
 $(1+\calL)^{-\beta}$, $\beta\!>\!0$, remain true.  The Ornstein-Uhlenbeck
 semigroup $T_t$ for $Y$-valued and real valued functions are related via
 \begin{equation} \Label{starstarstarK}
   T_{t}\clfY{f(\om)}{e}=\clfY{T_t f(\om)}{e}\hspace{1.5cm}
                        {\rm for\ all\ }f\!\in\!L^p(\mu;Y),\ e\!\in\!Y^*,
 \end{equation}
 In particular, if $a(\om)$ \marker{is representable in the sense
 of~(\ref{newstar})},\ then
 \[T_t a = \sum_i T_t (a, e_i)_Q\,e_i \]
 and similarly for $\calL$, etc.  Moreover, (i)---(iv) ,under obvious
 modifications, remain true.
%
 \subsection{Stochastic Differentiation} \Label{Stochdiff}
 Let
 \begin{equation} \Label{Sn}
  \calS_n=\mbox{\Large {$\{$}}\Phi=\vph(\delta l_1,\ldots,\delta l_n)\
         \mbox{\Large {$|$}}\ l_i\!\in\!W^*,\ i=1,\ldots,n,\ \
     \vph\!\in\!C^\infty_b(\Reals^n)\mbox{\Large {$\}$}}.
 \end{equation}
 For any $\Phi\!\in\!\calS_n$ represented as in~(\ref{Sn}), its gradient is
 the $W^*$--valued random variable
  \begin{equation} \Label{gradient}
  \grad{\Phi}=\sum_{i=1}^n\ \pd{\vph}{x_i}(\delta l_1,\ldots,\delta l_n)\,l_i,
 \end{equation}
 and this definition can be easily seen not to depend on $\Phi$'s particular
 representation. Denote
 \begin{equation}  \Label{S}
  \calS=\bigcup_{n=1}^\infty \calS_n.
 \end{equation}

 \noi
 The classical Sobolev completion of $\calS$ yields a space of functionals
 differentiable along $H$. In fact, other Sobolev spaces can be obtained
 by considering different subspaces of differentiability. Given a Banach space
 $(Z,\|\ \|_{_Z})$ continuously embedded in $W$\ \ (the elements of $Z$
 will be the directions of differentiability; cf. Remark~\ref{meaningful}
 below)\marker{$W^*\!\subset\!Z^*$\ and\ $\grad\Phi$ can be viewed as being\
 $Z^*$ valued; indeed
 \begin{equation} \Label{Zdual}
  \pair{Z}{z}{\grad{\Phi}(\om)}{Z^*}
    =\sum_{i=1}^n\ \pd{\vph}{x_i}(\delta l_1,\ldots,\delta l_n)\,
                                              \pair{Z}{z}{l_i}{Z^*}.
 \end{equation}
 Thus} for all $\Phi\in\calS$\ and\ $p\in[1,\infty)$ consider the
 Sobolev norms on\ $S$
 \begin{eqnarray}  \Label{Sobnorms}
  \|\Phi\|_{p,1;Z}&=&\left(\|\Phi\|^p_{L^p(\mu)} 
              +\|\grad{\Phi}\|^p_{L^p(\mu;Z^*)}\right)^{\frac{1}{p}}
 \end{eqnarray}
 and denote $\calS$'s completion according to this norm by \ $\Dp{Z}$.

 \noindent
 {\bf Example:}\ \ The Hermite polynomials\ \  $H_n(x)\!=\!\frac{(-1)^n}
    {\sqrt{n!}}\, e^{\frac{x^2}{2}}\mbox{\Large $\frac{d^n}{dx^n}$}
    \!\left(e^{-\frac{x^2}{2}}\right),\ \ \ n\!=0,1,2,\ldots$\ \ satisfy \ \
    $EH_n(X)H_m(X)=\delta_{n,m}$\ \ for\ $X\!\sim\!N(0,1)$\ \ and\ \
    $H_n'=\sqrt{n}H_{n-1}$\ \ for\ $n=1,2,\ldots$.
    Given an ONB\ ${e_n}$\ of $H$,\ and by Levy's criterion, the sequence
    of random variables\
    $a_m=\sum_{n=1}^m \frac{H_{2n}(\delta e_n)}{\sqrt{n}\log n}\in\calP_n$
    converges in\ $L^2$ and a.s. to, say, $a$\ and $\grad{a_m}=
    \sum_{n=1}^m \frac{\sqrt{2}}{\log n}\,H_{2n-1}(\delta e_n)e_n$.\ \
    This $H$--valued sequence does not converge in $L^2(\mu;H)$;
 If, however, for the case where the $W$ space is the completion of
 $H$ under the norm $\|u\|_{_W}= \left(\sum_i \left|\one \lip
 u,e_i\rip\right|^2\right)^\half$ or $\|u\|_W = \| Qu\|_{_H}$ where
 $Q$ is a Hilbert--Schmidt operator on $H$, $\grad a_m$ converges
 in $W$.  Therefore, in this case $a\in \DD_{2,1}^{W^*}$ but
 $a\not\in \DD_{2,1}^H$.  More generally, for any abstract Wiener
 space, $W$, we can embed a $W_0$ of the form defined above,\ i.e.\
 $H\subset W_0 \subset W$, and then $a\in \DD_{2,1}^{W_0^*}$.

 \begin{remark} \Label{meaningful} \rm \small
    \marker{In view of~(\ref{Zdual})\ }it is natural to think of $\grad{\Phi}$ being characterized by
    \[\pair{Z}{z}{\grad{\Phi}(\om)}{Z^*}
              =\lim_{\eps\to0}\frac{\Phi(\om\!+\!\eps z)-\Phi(\om)}{\eps}\]
    in some sense, however $\Phi(\om\!+\!\eps z)$ is meaningless unless
    $z\!\in\!H$ or $\Phi$ is sufficiently regular. For $Z\!\not\subset\!H$,
    thus, the space $\Dp{Z}$ consists of functionals which may be ``too regular"
    to be interesting in some applications. However, at the present they seem to
    be needed for the construction of flows on Wiener space as will be seen in
    Section~\ref{flowsection}.
  \end{remark}
  \vspace{.3cm}

 It is straightforward to verify that $\calS$ is dense in $L^p(\mu)$ and that
 the operator\ $\nabla$ is closeable with respect to \ $\|\ \|_{p,1;Z}$,\ with
 domain\ $\calS\!\subset\!L^p(\mu)$ and range in $L^p(\mu;Z)$.
 The space $\Dp{Z}$\ can thus be taken to be a dense subset of\ $L^p(\mu)$,
 and $\nabla$ has natural bounded linear extension\
 $\nabla^{\!^Z}:\Dp{Z}\longrightarrow L^p(\mu;Z^*)$.\
 When no confusion arises, $\nabla$'s \marker{superscript} may be omitted.
 In particular, if $\calP_n$ is the space of random variables obtained
 when $C^\infty_b(\Reals^n)$ is replaced in~(\ref{Sn}) by the family of
 polynomials in $n$ variables, and\ $\calP\!=\!\bigcup_n \calP_n$,\ then\
 $\calP\!\subset\!\Dp{W}$ and~(\ref{gradient}) still holds for any\
 $\Phi\!\in\!\calP$.

 \noindent
 If $(Z_1,\|\ \|_{Z_1})\imb(Z_2,\|\ \|_{Z_2})$\ then\
 $\|\Phi\|_{p,1;Z_1}\le\|\Phi\|_{p,1;Z_2}$\ and\ $\Dp{{Z_2}}$\ is continuously
 embedded in\ $\Dp{{Z_1}}$. In particular $\DD_{p,1}\!=\!\Dp{H}$\ is the
 classical Sobolev space in the Wiener context, and $\Dp{{W^*}}$ is a larger
 space consisting of Wiener functionals ``differentiable only along the $W^*$
 directions".

 \noi
 Finally, differentiation can also be defined for random variables taking
 values in a separable Banach space $Y$. Let\ $\calS(Y)$ (resp. $\calP(Y)$)\
 be $L^\infty(Y)$'s
 subset of elements having the form $\FFF=\sum_{k=1}^m\!\Phi_iy_i,$\ where
 $m\!\in\!\Nat$,\ \ $\Phi_i\!\in\calS$ (resp. $\calP$)\ \ and\ \
 $y_i\!\in\!Y,\ \ i\!=\!1,\ldots,\!m$.
 The gradient of such an\ $\FFF$ is defined to be
 \begin{equation*}
  \grad{\FFF}=\sum_{j=1}^m \grad{\Phi_i}\otimes y_i\ \in L^\infty(\mu;L(W,Y)).
 \end{equation*}
 (Here $\ W^*\!\otimes\!Y$\ is embedded naturally in\ $L(W,Y)$ by setting\
 $(l\otimes y)w\!=\!\lip w,l\rip\,y$). We then define, for a given\
 $W^*\imb Z\imb W$\ as above,\ $p\!\in\![1,\infty)$\ and\
 $\FFF\!\in\!\calS(Y)$,\ the Sobolev norms
 \begin{eqnarray*}
   \|\FFF\|_{p,1;Z}&=&  \left(\|\FFF\|^p_{L^p(\mu;Y)}
              +\|\grad{\FFF}\|^p_{L^p(\mu;L(Z,Y))}\right)^{\frac{1}{p}}
 \end{eqnarray*}
 and\ $\Dp{Z}(Y)$\ will be\ $\calS(Y)$'s completions according to these norms.
 The same monotonicity relations hold in the  differentiation space $Z$
 as in the scalar case, and $\nabla$ can be extended to a bounded operator\ \
 $\grad^{^Z}:\Dp{Z}(Y)\to L^p(\mu;L(Z,Y))$.

 \section{The Divergence Operator}  \Label{divsection}
\subsection{The Divergence of $\boldmath W$-valued r.v.'s} \Label{divsection1}

 \noi
 The standard definition of the divergence introduces it as an operator
 on suitable $H$--valued random variables (cf.~e.g.~\cite{M},
 \cite{U}). We now wish to extend it to $W$--valued
 random variables. In fact, with no extra effort, the same definition will
 also apply to $\Wss$--valued random variables, one advantage of
 which is pointed out in Remark~\ref{whyWstaraim}.
 \begin{definition} \Label{divdef}
  For $p\in[1,\infty)$ the space $\domdel{p}\!\subset\!L^p(\mu;W^{**})$ is
  defined to be the set of all $v\!\in\!L^p(\mu;\Wss)$\ for which there exists
  a random variable $\delta v\!\in\! L^p(\mu)$,\ the {\bf divergence} of\ $v$,\
  such that, for all\ $\Phi\!\in\!\calS$,
  \begin{equation}   \Label{scaduality}
   E\clfs{v}{\grad{\Phi}}=E\Phi\delta v.\ \ \
  \end{equation}
 \end{definition}
  In particular,
  \begin{equation}  \Label{Edelzero}
    E\delta v=0\hspace{1.5cm}\forall v\!\in\!\domdel{p}.
  \end{equation}
  Moreover, it follows form the $L^p$ duality theory (and \  $\calS$'s
  density in $L^p(\mu)$) that a necessary and sufficient condition
  for $v\!\in\!\domdel{p}$ is the existence of a finite positive
  constant $\gamma\!=\!\gamma(v)$ such that for all\ $\Phi\!\in\!\calS$
  \begin{equation} \Label{boundedlf}
    E\left|\clf{\grad{\Phi}}{v}\right|\le  \gamma\|\Phi\|_{L^q(\mu)}
  \end{equation}
  where $q$ is $p$'s conjugate exponent, $\frac{1}{p}\!+\!\frac{1}{q}\!=\!1$,
  in which case $\|\delta v\|_{L^p(\mu)}$ is the best possible
  constant\ $\gamma$\ in~(\ref{boundedlf}). In fact, if\ $p\!>\!1$\ and\
  $v\!\in\!\domdel{p}$,\ \ (\ref{scaduality}) and~(\ref{boundedlf}) will
  actually hold for all\ $\Phi\!\in\!\mathbb{D}_{q,1}^W$, in particular for all
  $\Phi\!\in\!\calP$.

 \begin{remarks} \Label{rmks}\rm \small a) The operator\ $\delta$\ and its domain\ $\domdel{p}$\
   are classical objects in the context of\ $H$--valued random variables.
   Obviously, if in our setup\ $v$ happens to take its values in $H$ a.s.
   (in which case the pairing in~(\ref{scaduality}) and in~(\ref{boundedlf})
   become\ $(v,\grad{\Phi})_H$\,) the definition of\ $\delta$ reduces to the
   classical one.

  \noi b) There are no new deterministic elements in\ $\domdel{1}$. Indeed,
  assume that a (deterministic)\ $v\!\in\!\Wss$\ belongs to\
  $\domdel{1}$, that is, satisfies~(\ref{boundedlf}) for all
  $\Phi\!\in\!\calS$\ and\ $p\!=\!1$. For each\ $l\!\in\!W^*$\ denote\
  $\Phi_l\!=\!\varphi(\delta l)\in\calS$,\ for some fixed\
  $\varphi\!\in\!C_b^{\infty}(\mathbb{R})$, strictly increasing,
  (e.g.\ $\varphi(x)\!=\!\arctan(x)$)\ so that \ $a:=E\varphi'(Z)\!>\!0$,\
  where\ $Z\!\sim\!N(0,1)$. For these test functions,\
  $\pair{W^*}{\grad\Phi_l}{v}{\Wss}\!=\!\varphi'(\delta l)\clfs{v}{l}$,\
  and it thus follows from~(\ref{boundedlf}) that
  \begin{equation}
   \sup_{l\in W^*,\ |l|_H=1} \left|\clfs{v}{l}\right|
                      \le\frac{\gamma\,\|\varphi\|_{\infty}}{a}<\infty.
  \end{equation}
  This implies that\ $v$ can be extended as a bouned linear functional on\ $H$,\
  i.e.\ $v\!\in\!H$. In other words, there are no deterministic elements in $W$,
  \ or even of\ $\Wss$, that possess a divergence without being in $H$, and thus
  already having a divergence in the classical sense.
 \end{remarks}

 In view of Remark~\ref{rmks}b) one might wonder whether there are
 any non $H$--valued elements in\ $\domdel{1}$\ at all. The
 following example answers this question affirmatively.

 \noi
 {\bf Example:\ } For a given ONB\ $\{e_i\}$\ in\ $H$, let
   $$v=(\delta e_2)e_1-(\delta e_1)e_2+(\delta e_4)e_3-(\delta e_3)e_4
    \pm\ldots\ \ \ .$$
   By the Ito--Nisio theorem,\ $v(\omega)$\ is a measure preserving
   transformation of\ $\omega$\ and thus isn't supported on $H$.
   We claim that\ $v$\ possesses a divergence, which moreover is a.s.\ $0$.
   This follows from the obvious fact that \
   $\delta\left[\,(\delta e_{2k})e_{2k-1}-(\delta e_{2k-1})e_{2k}\,\right]=0$\
   for all $k\!\in\!\mathbb{N}$, and from Lemma~\ref{approx1}
   below which extends this equality to the infinite sum.

   \noi
   It is interesting to note that, on the other hand,\
   $\tilde{v}(\omega)=\omega=\sum_{i=1}^\infty \delta(e_i)\,e_i$,
   which is (trivially) a measure preserving transformation of $\omega$
   as well, does {\it not} have a divergence.

%

   \begin{lemma}  \Label{av}
   Let $\alpha\!\in\!\DD_{p_1,1}^{W}$\, \ $v\!\in\!\dom_{p_2}\delta$\ and\
   $\frac{1}{p}\!=\!\frac{1}{p_1}\!+\frac{1}{p_2}$.\ Then\
   $\alpha v\!\in\!\dom_p\delta$\ and\
   \begin{equation}  \Label{deltaav}
    \delta(\alpha v)=\alpha\delta v-\clfs{v}{\grad{\alpha}}\ .
   \end{equation}
  \end{lemma}
  (As observed in Remark~\ref{meaningful}, the family of $\alpha$'s allowed
  in this lemma is only slightly more general than $\calS$.\ This result will
  be applied in Lemma~\ref{secmom} at the end of this section).

  \noi
  \pf
     For every $\Phi\!\in\!\calS$
     \begin{eqnarray*}
      E\clfs{\alpha v}{\grad{\Phi}}&=&E\clfs{v}{\alpha\grad{\Phi}}
              = E\clfs{v}{\grad{(\alpha\Phi)}}-E\clfs{v}{\Phi\grad{\alpha}}\\
             &=&E(\alpha\Phi\delta v)-E\Phi\clfs{v}{\grad{\alpha}}
              = E\Phi\left(\alpha\delta v-\clfs{v}{\grad{\alpha}}\right)
     \end{eqnarray*}
     which proves the result.
   \qed
 \begin{lemma}  \Label{approx1}
  Let $p\!\in\![1,\infty)$\ \ and\ \
  $\{v_n\}_{n=1}^\infty\!\subset\!\dom_p\delta.$\ \ If

  \hspace*{.5cm}\begin{minipage}[t]{16cm}
  \begin{itemize}
   \item[i)] $v_n\longrightarrow\hspace{-.77cm}
         \raisebox{-.14cm}{\mbox{\tiny $n\!\!\to\!\!\infty$}}\hspace{.15cm}v$
    \ \ weakly in $L^p(\mu;W)$\ \ \ \ and
   \item[ii)] $\{\delta v_n\}_{n=1}^\infty$\ \ is bounded in $L^p(\mu)$, \ \
    i.e. $\exists M\!<\!\infty$\ \ such that\ \
    $\|\delta v_n\|_{L^p(\mu)}\!\le\!M$\ \ for all\ $n$,
  \end{itemize}
  \end{minipage}

  \noi then\ \ $v\!\in\!\dom_p\delta$\ \ \ \ and\ \ \
  $\delta v_n\longrightarrow\hspace{-.77cm}
   \raisebox{-.14cm}{\mbox{\tiny $n\!\!\to\!\!\infty$}}
   \hspace{.2cm}\delta v$ \ \ weakly in\ $L^p(\mu)$.\ \ \
   In particular\ $\|\delta v\|_{_{L^p(\mu)}}\le M$.
 \end{lemma}
 \pf Clearly\ $v\!\in\!L^p(\mu;W)$.\
     Defining $q$ as usual by\ $\frac{1}{p}\!+\!\frac{1}{q}\!=\!1$,
     \[ \left|E\clf{\grad{\Phi}}{v_n}\right|
                \le\|\delta v_n\|_{L^p(\mu)}\|\Phi\|_{L^q(\mu)}
                \le M\|\Phi\|_{L^q(\mu)}\hspace{2cm}
                \forall n\!\in\!\Nat,\ \ \ \ \forall \Phi\in\calS,\]
     so that by~($i$)\ \ $v$ satisfies~(\ref{boundedlf})\ and thus\
     $v\!\in\dom_p\delta$. Moreover, for every\ $\Phi\!\in\!\calS$,
     \begin{equation} \Label{wlimit}
       E\Phi\delta v=E\clf{\grad{\Phi}}{v}
                  =\lim_{n\to\infty}E\clf{\grad{\Phi}}{v_n}
                  =\lim_{n\to\infty}E\Phi\delta v_n.
     \end{equation}
     Since\ $\calS$\ is dense in $L^q(\mu)$\ and\
     $\{\delta v_n\}_{n=1}^\infty$\ is\ bounded in $L^p(\mu)$,\
     the end terms of~(\ref{wlimit}) are equal for all $\Phi\!\in\!L^q(\mu)$,\ \
     in other words\ $\delta v_n\longrightarrow\hspace{-.77cm}
     \raisebox{-.14cm}{\mbox{\tiny $n\!\!\to\!\!\infty$}}
     \hspace{.2cm}\delta v$ \ \ weakly in\ $L^p(\mu)$.
 \qed\\
 \begin{corollary}  \Label{DDindomdel}
  \ \ $\Dp{W}(\Wss)\subset\dom_p\delta$\ \ \ \ and\ \
  $\delta:\Dp{W}(\Wss)\longrightarrow L^p(\mu)$\ \ is a bounded linear
     \marker{operator}.
 \end{corollary}
 \pf
  Recall  that $\delta:\Dp{H}(H)\to\Lp$
  is a bounded linear operator, i.e. there exists a finite constant\
  $C$\ such that\ $\|\delta \Phi\|_{\Lp}\!\le\!C|\Phi|_{p,1}$\ for every
  $\Phi\!\in\!\Dp{H}(H)$.

  \noi
  Let $\FFF\in\Dp{W}(\Wss)$. By definition there exists a sequence\
  $(\FFF_{_n})_{_n}\!\subset\!S(\Wss)$\ such that\ $\FFF_{_n}\to \FFF$\
  in\ $\Dp{W}(\Wss)$.\
  With no loss of generality we may in fact assume for each $n$ that\
  $\|\FFF_n\|_{p,1}\!\le\!\|\FFF\|_{p,1}$\ and (since\ $S(H)$\ is dense
  in\ $S(\Wss)$\ in the\ $\|\cdot\|_{p,1}$\ norm) that\
  $\FFF_{_n}\!\in\!S(H)$.\ Thus
  \[ \|\delta \FFF_{_n}\|_{_{\Lp}}\le C|\FFF_{_n}|_{_{p,1}}
      \le C\|\FFF_{_n}\|_{_{p,1}}\le C\,\|\FFF\|_{_{p,1}}\hspace{1cm}\forall n\ ,\]
  so that by Lemma~\ref{approx1}\ (obviously\ $\FFF_{_n}\to \FFF$\ weakly
  in\ $\Lp$) we conclude that\ $\FFF\!\in\!\dom_p\delta$\ and that
  $\|\delta \FFF\|_{_{L^p(\mu)}}\!\le\!C\|\FFF\|_{_{p,1}}$.
 \qed
 \bigskip

  \noi
  In some cases we have the following convenient approximation.
  Given a smooth ONB\ \ $\calE\!=\!\{e_i\}_{i=1}^\infty$ of $H$\,
  assume $v\!\in\!L^p(\mu;W)$ has the representation
  \begin{equation} \Label{series}
   v=\sum_{i=1}^\infty \alpha_ie_i
  \end{equation}
  that is,\ \ $v\!=\!\lim_{_{\!n\to\infty}}\!v(n)$\ \ in
\ $\LpW$,\  where\
  $v(n)\!=\!\sum_{i=1}^n\alpha_ie_i$.\ Note that\
  $\alpha_i\!=\!\clf{e_i}{v}\!\in\!\Lp$\ for each $i$,\ so that it is
  possible, for each\ $n\!\in\!\Nat$,\ to define the\ $H$--valued projections
  \begin{equation} \Label{projection}
   v_n=v_n^{\calE}
      =\sum_{i=1}^n\EFn{\alpha_i}e_i=\EFn{v(n)}
  \end{equation}
  where $\calF_n\!=\!\calF_n^{\calE}\!
    =\!\sigma(\delta e_1,\ldots,\delta e_n)$,\ the sigma--algebra generated
  by $\delta e_1,\ldots,\delta  e_n$, and the conditional expectation on
  the right is for $W$--valued random variables.
  \begin{proposition}  \Label{martingale}
   Let\ $p\!\in\![1,\infty)$ and assume that\  $v\!\in\!\dom_p\delta$
   is represented as in~(\ref{series}). Then \vspace{-.5cm}

   \begin{itemize}
    \item[i)]   $\lim_{_{n\to\infty}}\!v_n=v$\ \ in\ $\LpW$
    \item[ii)]  $v_n\!\in\!\dom_p\delta$\ \ and\ \ $\delta v_n\!
                =\!\EFn{\delta v}$
    \item[iii)] $\lim_{_{n\to\infty}}\delta v_n=\delta v$\ \
                                            a.s. and in\ $\Lp$.
   \end{itemize}
  \end{proposition}
  \pf
   The first claim follows from
   \begin{eqnarray}
    \left\|v-v_n\right\|_{_{_{\LpW}}}
           &\le& \left\|v-\EFn{v}\right\|_{_{_{\LpW}}}
                +\left\|\EFn{v-v(n)}\right\|_{_{_{\LpW}}}\nonumber\\
           &\le& \left\|v-\EFn{v}\right\|_{_{_{\LpW}}}
                +\left\|v-v(n)\right\|_{_{_{\LpW}}}.\Label{Lpconv}
   \end{eqnarray}
   The first term converges to zero by the ($W$--valued) martingale
   $L^p$ convergence theorem, while the second term does so by
   assumption.

   \noi
   Next, note that\ $\EFn{\Phi}\!\in\!\calS$\ for an arbitrary $\Phi\!\in\!\calS$\
   \ and\ $\grad{\EFn{\Phi}}\!=\!\conde{\grad{\Phi}}{\calF_n}$.
   Then
   \begin{eqnarray*}
    E\clf{\grad{\Phi}}{v_n}&=&E\clf{\EFn{\grad{\Phi}}}{v(n)}
                       = E\clf{\grad{\EFn{\Phi}}}{v(n)}\\
                      &=&E\clf{\grad{\EFn{\Phi}}}{v}
                       = E\left(\delta v\,\EFn{\Phi}\right)
                        =E\left(\EFn{\delta v}\Phi\right)
   \end{eqnarray*}
   which proves (ii).

   \noi
   Finally (iii) follows from (ii) by the martingale convergence theorem
   applied to\ \ $\delta v_n\!=\!\EFn{\delta v}$\, in both the a.s. and
   $L^p$ senses.
  \qed

  \noindent
  In Section~\ref{flowsection} we shall need the following extension of
  the Proposition's first statement to $v$'s parametrized by some positive
  measure space $(I,{\cal I},\lambda)$.
  \begin{corollary} \Label{martingalecor}
   For a given \ $p\!\in\![1,\infty)$\ \ assume that\ \
   $v_t(\om)\!\in\!L^p(I\!\times\!W,\lambda\!\times\!\mu\,;\,W)$\ \ satisfies
   \[ \lim_{n\to\infty}\int_I
        E\left\|v_t-\sum_{i=1}^n\alpha_{i_t}e_i\right\|^p_Wd\lambda(t)=0\]
   for some smooth ONB\ \ $\calE\!=\!\{e_i\}_{i=1}^\infty$\ \ of $H$
   (with $\alpha_{i_t}=\clf{v_t}{e_i}$).\  Then
   \[ \lim_{n\to\infty}
         \int_I
         E\left\|v_t-\sum_{i=1}^n\EFn{\alpha_{i_t}}e_i\right\|^p_W
                                                        d\lambda(t)=0.\]
  \end{corollary}
   Its proof essentially repeats the one of Proposition~\ref{martingale}$i)$
   except that one must add a dominated convergence argument (to the integral
   over $I$) for the first term in~(\ref{Lpconv})'s appropriate extension to
   converge to zero.
  \vspace{.2cm}

  \begin{definition}  \Label{dfreexact}
   Let $p\!\in\![1,\infty)$.
   \vspace{-.5cm}

   \begin{itemize}
    \item[a)] An element\ \ $v\!\in\!\domdel{p}$\ \ is said to be
             {\em divergence-free} if\ \ $\delta v\!=\!0$,\ \ i.e.
             if\ \ $E\clfs{v}{\grad{\Phi}}\!=\!0$\ \ for all\
                                                $\Phi\!\in\!\calS$.\
             The class of all such divergence-free\  is denoted by\
             $\domo{p}$.
    \item[b)] An element\ \ $u\!\in\!L^p(\mu;H) \!\subset\!\LpW$ \ \ is said
     to be {\em exact} 
     if there exists a\ \ $\Psi\!\in\!\Dp{H}$\ \  such that\ \ $u=\grad{\Psi}$.\ \
     This class of exact\ $H$--valued random variables is denoted\ \ $\Le{p}$.
   \end{itemize}
  \end{definition}
   \begin{lemma}  \Label{domolemma}
    Let $p\!\in\![1,\infty)$.\vspace{-.3cm}

    \begin{itemize}
     \item[a)] If\ $u\!\in\!\domo{p}$\ \ \ then\ \
               $(1\!+\!\calL)^{^{-\beta}}\!u\in\domo{p}$\ \ \
               for every $\beta\!\ge\!0$.
     \item[b)] \ $\domo{p}\cap\Le{p}=\{0\}$
    \end{itemize}
   \end{lemma}
   \pf
    {\bf a)} \
 By \req{starstarK}\ \ $(1\!+\!\calL)^{-\beta}u \in L^p(\mu;W^{**})$.\ \ \
 By~\req{starstarK},~\req{starstarstarK} and $T_t$'s self-adjointness,
 \ \ for any $\Phi\!=\!\varphi(\delta e_1,\dotsc,\delta e_n)\in\calS$,\ \ \
 $\nabla\Phi\!=\!\sum_1^n \varphi_i e_i\in\calS(W^*)$,\ \ \ (here\ $\varphi_i\!
   =\!\frac{\partial \varphi}{\partial x_i}(\delta e_1,\dotsc,\delta e_n)$),
   \ \ and
 \begin{align*}
 E\clfss{\grad\Phi}{(1+\calL)^{-\beta}u}
 & = \frac{1}{\Gamma(\beta)} \: E \int_0^\infty t^{\beta-1} e^{-t}
 \clfsss{T_tu}{\sum_1^n \varphi_i e_i} \, dt \\
 & = \frac{1}{\Gamma(\beta)} \: E \int_0^\infty t^{\beta-1} e^{-t}
 \sum_1^n \Bigl( \clfsss{T_tu}{e_i} \,\varphi_i \Bigr)\, dt\\
 & =\frac{1}{\Gamma(\beta)} \: E \int_0^\infty t^{\beta-1} e^{-t}
 \sum_1^n \Bigl(T_t \lip u, e_i \rip \Bigr)\, \varphi_i \, dt\\
 & =\frac{1}{\Gamma(\beta)} \:  \int_0^\infty t^{\beta-1} e^{-t}
 \sum_1^n  E \Bigl(T_t \varphi_i \lip u, e_i \rip \Bigr)\, dt\\
 & = E \clfs{u}{(1+\calL)^{-\beta} \grad \Phi}\\
 & = E\clfs{u}{\nabla \calL^{-\beta} \Phi}\hspace{1.5cm}
                     \text{(by (iv) in subsection~\ref{Stochdiff})} \\
 &= 0\,.
 \end{align*}

    {\bf b)} If $u\!=\!\grad{\Psi}$\ and\ $\delta u\!=\!0$, then
     $0\!=\!\delta\grad{\Psi}\!=\!\calL \Psi$. Thus 
     \ $\Psi\!=\!E\Psi$ a.s.,\ \ so that \ $u\!=\!0$.
    \qed

   \noi
   The following proposition essentially states that the only ``new"\
   $W$--valued vector fields with divergence are divergence free.
  \begin{proposition} \Label{vfdecomp}
   Let\ $p\!\in\![1,\infty)$. Each\ $v\!\in\dom_p\delta$\ can be uniquely
   decomposed as a sum\ \ $v\!=\!v^0\!+\!v_{\rm e}$\ where\ $v^0$\ is divergence
   free and\ $v_{\rm e}$\ is exact (with divergence). Equivalently
   \[ \dom_p\delta=\dom^0\delta\oplus \left(\Le{p}\cap\dom_p\delta\right). \]
  \end{proposition}
  \pf
   Let\ $v\!\in\dom_p\delta$. By~(\ref{Edelzero}) and the remarks following
   \req{starstarK}, $v_e\!=\!\grad{\left(\calL^{-1}\delta v\right)}$\ \ is a
   well defined element of\ $\Le{p}\cap\dom_p\delta$. In order to prove the
   decomposition, we need to check that\ $v^0\!:=\!(v-v_e)\!\in\!\domo{p}$.\
   Indeed, as remarked above, $v_e\!\in\!\dom_p\delta$.\ Moreover, by\ (ii) of
   subsection~(\ref{calL}),
   \[ \delta v^0=\delta v-\delta\grad{\left(\calL^{-1}\delta v\right)}=0.\]
   The uniqueness follows directly from Lemma~\ref{domolemma}b.
  \qed
 \begin{remark}  \rm \small \Label{Euclidean}
  Heuristically, vector fields which possess divergence generate flows. This will
  be formalized, under appropriate assumptions, in Section~\ref{flowsection}
  noting in addition that the flows generated by divergence free vector fields
  are measure preserving (``rotations") while those generated by $H$--valued
  vector fields (``shifts") have been already studied, for example
  in~\cite{Cru},~\cite{P}\ and~\cite{UZ99}. What Proposition~\ref{vfdecomp}
  suggests is that a general vector field which generates flows can be decomposed
  into a ``rotation" generating component and a ``shift" generating  component.
 \end{remark}

 \marker{The formula for the classical divergence's second moment has its
   counterpart for $W$--valued variables as well, but since it involves operators
   and their divergence, it is deferred till the next subsection
   (Lemma~\ref{secmom}).}
 \subsection{The Divergence of Operator Valued r.v.'s}  \Label{divsection2}

   In the next section divergence--free vector fields will be
   characterized as the divergence of an antisymmetric operator.
   This is, at this stage, only a formal declaration. The remainder of
   this section is dedicated, therefore, to precise what is meant by
   the divergence of an operator and to present some of its properties.

  \noi
  Indeed, the classical divergence is defined for random variables taking
  values not only in $H$ but in $H$'s tensor powers as well. Thus, for
  example, the divergence $\ddelta \AAA$ of a random Hilbert Schmidt
  operator $\AAA(\om)$ in $H$ is characterized by
  \begin{equation} \Label{HSdiv}
   E(\AAA,\grad{\FFF})_{_{H^{\!\mbox{\tiny $\otimes 2$}}}}
       =E(\ddelta \AAA,\FFF)_{_H}\hspace{1cm}\forall \FFF\!\in\!\calS(H).
  \end{equation}
  Here $(\AAA,\BBB)_{_{H^{\!\mbox{\tiny $\otimes 2$}}}}={\rm tr \AAA\BBB^T}
   =\sum_{i=1}^\infty (\AAA e_i,\BBB e_i)_{_H}$ (for any ONB $\{e_i\}$)\
  is the natural  inner product of two Hilbert Schmidt operators on $H$.

  \noi
  The random operators we wish to generalize the divergence\ $\ddelta$\ to
  will have $W^*$ as their domain and a fixed arbitrary Banach space as their
  range (instead of\ $\mathbb{R}$\ as in subsection~\ref{divsection1}).
  To carry out this generalization, recall first the definition of the
  trace ${\rm tr}\,T$\ of\ an operator\ $\TTT\!\in\!L(W^*,\Wss)$,\ namely\
  $\sum_{i=1}^\infty \clfs{\TTT e_i}{e_i}$ if this sum exists for every smooth
  ONB $(e_i)$\ of\ $H$,\ and is the same for all such bases.
  In particular, every finite rank operator from\ $W^*\ to\ \Wss$  has a trace.\\

  \noi
  Henceforth, $Y$ will be a fixed Banach space,\ $p\!\in\![1,\infty)$\, and\
  $\frac{1}{p}\!+\!\frac{1}{q}\!=\!1$.
  \begin{definition} \Label{Domp}
  \hspace*{.35cm}Define\ \ $\Domdel{p,Y}\!=\!\Domdel{p}\!$\ \ to be the set of all\ \
  $\KKK\!\in\!L^p(\mu;L(W^*,Y))$\ \ for which there exists a\ \
  $\ddelta \KKK\!\in\!L^p(\mu;Y)$,\ \ the {\bf divergence} of\ $\KKK$,\ \
  such that
  \begin{equation} \Label{opdual}
   E\,{\rm tr}\left(\KKK^T\grad^{W^*}{\!\FFF}\right)
         =E\pair{Y}{\delta \KKK}{F}{Y^*}
  \end{equation}
  for all\ $\FFF\!\in\!\calS(Y^*)$.
  \end{definition}
  \begin{remark} \ \ \Label{finiterank}  \rm \small
  In~\req{opdual}, $\grad^{\!^{W^*}}\!\!\FFF$\ 
  has finite rank a.s., so that its left hand side makes sense.
  If\ $\KKK$\ itself has a deterministic finite dimensional range,~\req{opdual}
  will also hold for all\ $F\!\in\!\mathbb{D}_{q,1}^{W^*}(Y^*)$.\ Indeed, both
  of its terms pass to the limit when\ $F_n\!\to\!F$\ 
  \ \ \ ($F_n\!\in\!\calS(Y^*)$).
  \end{remark}
  \vspace{.3cm}

  The uniqueness of\ $\ddelta\KKK$\ is a consequence of\ $S(Y^*)$'s density
  in\ $L^q(\mu;Y^*)$.\ Moreover, if\ $\KKK\!\in\!\Domdel{p,Y}$\ and\ $Y$\ is
  continously embedded in another Banach space $Y_1$, then~\req{opdual}
  obviously holds for all\ $F\!\in\!\calS(Y_1^*)\!\subset\!\calS(Y^*)$. Thus\
  $\Domdel{p,Y}\!\subset\!\Domdel{p,Y_1}$\ and\ $\ddelta \KKK$\ is the same\
  $Y$--valued random variable whether\ $\KKK$'s range is taken to be\ $Y$\ or\
  $Y_1$. That is the reason why it isn't necessary to include the subscript\
  $Y$ in the notation of\ $\ddelta$.

  Just as in the scalar case a necessary and sufficient condition for
  $\KKK\!\in\!\Domdel{p,Y}$ is the existence of a finite positive constant
  $\gamma\!=\!\gamma(\KKK)$ such that for all\ $\FFF\!\in\!\calS(Y^*)$
  \begin{equation} \Label{boundedlf1}
   |E\,{\rm tr}\left(\KKK^T\grad^{W^*}\FFF\right)|
                           \le  \gamma\|\FFF\|_{L^q(\mu;Y^*)}
  \end{equation}
  in which case $\|\ddelta \KKK\|_{L^p(\mu;Y)}$ is the best possible
  constant\ $\gamma$\ in~(\ref{boundedlf1}).

  We denote\ $\Domdel{p}=\Domdel{_{p,\Wss\,}}$. Indeed, in Sections~\ref{domo}
  and~\ref{flowsection}, $Y$\ will typically be\ $\Wss$\ and together with\
  $\delta \KKK$ we shall need to consider\ $\delta \KKK^T$\ as well.
  (By a slight abuse of notation,\ $\KKK^T$\ actually stands for\
  $\KKK^T\left|^{}_{_{W^*}}\right.$). Recall that\ $L(W^*,\Wss)$\ can be also
  seen as the space $M_2(W^*)$ of bilinear forms in\ $W^*$,\ and in this
  interpretation, $\KKK^T(l_1,l_2)\!=\!\KKK(l_2,l_1)$). In particular, then,
  $\ddelta$'s domain contains $W\!\otimes\!W$ in this case, as stated in the
  abstract.\\


  \noi
  A useful connection between this divergence\ $\ddelta$\ and its scalar
  counterpart\ $\delta$\ is the following
  \begin{lemma}\Label{weakversion}
    An element\ $\KKK\!\in\!L^p(\mu;L(W^*,Y^{**})$\ belongs to
    $\Domdel{p,Y^{**}\,}$\ \ if and only if \ \  $\KKK^Tl\in\domdel{p}$\ \
    for every $l\!\in\!Y^*\!\subset\!Y^{**}$\ and for some\ $C\!>\!0$
    \begin{equation} \Label{Cbound}
     \|\ddelta \left(\KKK^Tl\right)\|_{_{L^p(\mu)}}\!\le\!C\|l\|_{_{Y^*}}
                                        \hspace{1.5cm}\forall\ l\!\in\!Y^*.
    \end{equation}
    In this case
   \begin{equation} \Label{weak}
    \delta(\KKK^Tl)=\pair{Y^*}{l}{\ddelta \KKK}{Y^{**}}\hspace{1.5cm}{\rm a.s.}
   \end{equation}
  and more generally, for any\ $F\!\in\!\calS(Y^*),\ \ \ \
  \KKK^TF\!\in\domdel{p}$\ \ \ and
 \begin{equation} \label{weakb}
  \delta(\KKK^TF)
      =\pair{Y^*}{F}{\ddelta \KKK}{Y^{**}}-{\rm tr}(\KKK^T\grad^{^{\!\!Y^*}}\!F)
 \end{equation}
  \end{lemma}
  \pf Throughout the proof, $\grad$\ will stand for\ $\grad^{^{\!\!Y^*}}$.\ \ \
      First assume that $\ddelta \KKK$ exists. For any $\Phi\!\in\!\calS$,\
      $l\!\in\!Y^*$ and denoting $G\!=\!\Phi l\!\in\!\calS(Y^*)$, it is
      straightforward to verify that
   \begin{equation} \Label{tracetp}
      {\rm tr}\,(\KKK^T\grad{G})=\clfs{\KKK^Tl}{\grad{\Phi}}\,.
   \end{equation}
   Then
   \begin{eqnarray*}
    E\clfs{\KKK^Tl}{\grad{\Phi}}=E\,{\rm tr}\,(\KKK^T\grad{G})
        =E\pair{Y^*}{G}{\ddelta \KKK}{Y^{**}}
        =E\Phi\pair{Y^*}{l}{\ddelta \KKK}{Y^{**}}
   \end{eqnarray*}
   and $\Phi\!\in\calS$ being arbitrary, $\ddelta(\KKK^Tl)$\ exists
   and~(\ref{weak}) holds, from which~(\ref{Cbound}) follows directly.

   \noindent
 %
   In the converse direction, it follows from~(\ref{Cbound}) that there exists
   a $\Delta_{_\KKK}\!\in\!L(Y^*,L^p(\mu))\approx L^p(\mu;Y^{**})$\ (we shall
   indeed relate to\ $\Delta_{\KKK}$\ as an element of\ $L^p(\mu;Y^{**})$)\
   such that for all $l\!\in\!Y^*$
   \begin{equation} \Label{Delta}
    \ddelta(\KKK^Tl)=\pair{Y^*}{l}{\Delta_{_\KKK}}{Y^{**}}\hspace{1cm}{\rm a.s.},
   \end{equation}
   so that, for any\ $\FFF=\sum_{j=1}^m \Phi_jl_j\!\in\!\calS(Y^*)$\
   \begin{eqnarray*}
    E{\rm tr}\,(\KKK^T\grad{\FFF})
      &=&\sum_{j=1}^m {\rm tr}\,\KKK^T\grad(\Phi_j\,l_j)
         =\sum_{j=1}^m E\clfs{\KKK^Tl_j}{\grad{\Phi_j}}\\
      &=&\sum_{j=1}^m E\delta(\KKK^Tl_j)\Phi_j
         =E\sum_{j=1}^m \pair{Y^*}{l_j}{\Delta_{_\KKK}}{Y^{**}}\Phi_j\\
      &=&\Phi\pair{Y^{*}}{\FFF}{\Delta_{_\KKK}}{Y^{**}}.
   \end{eqnarray*}
  Thus\ $\ddelta \KKK$\ exists by definition and is actually\
  $\Delta_{_\KKK}$, (\ref{Delta}) being nothing else but~(\ref{weak}).

  \noindent
  Turning to \req{weakb}, and with\ \ $F\!=\!\Phi\,l$,\ \ ($\Phi\!\in\!\calS$\
  and\ $l\!\in\! Y^*$),\ \ it follows from Lemma~\ref{av} that
  \[ \delta(\KKK^TF)\ =\ \delta(\Phi\KKK^Tl)
                    \ =\ \Phi\,\delta(\KKK^Tl)-\!\clfs{\KKK^Tl}{\grad{\Phi}}\]
  which proves the claim in view of~(\ref{weak})\ and~(\ref{tracetp}).
  \qed
 \begin{remark} \Label{whyWstaraim} \rm \small
  This lemma might shed some light on two questions concerning the r\^ole of\
  $Y^{**}$\ in general and\ $\Wss$\ in particular. First, the left hand side
  of~(\ref{weakversion}) requires the (scalar) divergence\ $\delta$ to be
  defined for random variables taking values in\ $\Wss$,\ not only in\ $W$.\
  Secondly, the vector divergence\ $\ddelta \KKK$\ of\ $\KKK$\ (initially
  identified in the proof as\ $\Delta_{_{\KKK}}$)\ must be allowed to take values
  in\ $Y^{**}$,\ not only in\ $Y$,\ for the lemma to hold.
 \end{remark}
%
%
  \begin{corollary}  \Label{extdelprops}
   For any $p\!\ge\!1$
   \mbox{}\\ \vspace{-1cm}

   \begin{itemize}
    \item[a)] $\Dp{W}(L(W^*,\Wss))\!\subset\!\ddomdel{p}$.
    \item[b)]
     If\ $v\!\in\!\domdel{p}$\ and\ $y\!\in\!\Wss$\ then\ \
     $v\otimes y\!\in\ddomdel{p}$\ \ and\ \ $\ddelta(v\otimes y)=\delta(v)y$.
    \item[c)] Let\ $\alpha\!\in\!\DD_{p_1,1}^{W}$,\ $\AAA\!\in\!\ddomdel{p_2}$\ \
              and\ $\frac{1}{p}\!=\!\frac{1}{p_1}\!+\frac{1}{p_2}$.\ \
    Then\ $\alpha \AAA\!\in\!\domdel{p}$\ \ and\ \
    \begin{equation}  \Label{cA}
     \ddelta \alpha\AAA=\alpha\ddelta \AAA-\AAA\grad^{\!^{W^*}}\alpha
             \hspace{1.5cm}{\rm a.s.}
    \end{equation}
   \end{itemize}
  \end{corollary}
%
%
 \pf Let $l\!\in\!W^*$.\ If\ $\KKK\!\in\!\Dp{W}(L(W^*,\Wss))$\ it is
  straightforward to verify that the same is true for $\KKK^T$, and therefore
  that\ $\KKK^Tl\!\in\!\Dp{W}(\Wss)$,\ so that\ $\KKK l\!\in\!\domdel{p}$ by
  Corollary~\ref{DDindomdel}. Since $l\!\in\!W^*$ was arbitrary, a) follows
  from Lemma~\ref{weakversion}. As for b), for any $l\!\in\!W^*$
  \begin{eqnarray*}
   \clfs{\ddelta v\otimes y}{l}
        & = &\delta\left((v\otimes y)^T l\right)
        \ =\ \delta(\clfs{y}{l}v)\\
        & = &\clfs{y}{l}\delta v\ =\ \clfs{(\delta vy)}{l}
  \end{eqnarray*}
  which proves the claim.

  \noindent
  Finally for c), assume without loss of generality that $\alpha\!\in\!\calS$.
  and let\ $l\!\in\!W^*$. Denoting\ $F=\alpha\,l$, it follows from~(\ref{weakb})
  that
   \[ \clfs{\delta(\alpha\AAA)}{l}=\delta(\AAA^TF)
           =\clfs{\delta\AAA}{\alpha l}-\clfs{\AAA\grad{\alpha}}{l} \]
     from which~(\ref{cA}) follows, again since $l\!\in\!W^*$ was arbitrary.
 \qed
 \bigskip

 \noi
  We conclude this section with an extension of a classical second
  moment identity.
    \begin{lemma}   \Label{secmom}
   Let $G\!\in\calS(W^*)$\ and\ $u\!\in\!\mathbb{D}^{^{W^*}}_{1,1}(\Wss)$.\
   Then
    \begin{equation} \Label{dudG1}
      E(\delta u\delta G)
         =E\clfs{u}{G}\!+\!E{\rm tr}\left(\grad G\,\grad^{^{\!\!W^*}}\!u\right).
    \end{equation}
     If, moreover,\ $\grad u\!\in\!\Domdel{1}$\ \ it will also hold that
    \begin{equation} \Label{dudG2}
          E(\delta u\delta G)
         =E\clfs{u}{G}\!+\!E\clfs{\delta(\grad u)^T}{G}
    \end{equation}
  \end{lemma}
   \pf
    Clearly\ $\alpha:=\delta G$\ belongs to $\calP$,\ so that by Lemma~\ref{av}
    \begin{eqnarray*}
     E(\delta u \delta G)
       &=&E\delta(\delta G\,u)+E\clfs{u}{\grad{\delta G}}\\
       &=&\ \ 0\ \ + \ E\clfs{u}{G}+E\clfs{u}{\delta((\grad{G})^T)}.
    \end{eqnarray*}
    In the \marker{last} equality we have used the well known identity\ \
    $\grad{\delta G}\!=\!G+\delta((\grad{G})^T)$ (in which the
    third term is, in fact, $W^*$--valued).
    The identity~\req{dudG1} will then follow by applying~\req{opdual} to the
    last term, together with Remark~\ref{finiterank}\ (here,\ $Y\!=\!W^*$).
    A second application of~\req{opdual} yields~\req{dudG2}, this time
    with\ $Y\!=\!\Wss$\ and thinking of\ $G$\ as an element of\
    $\calS(W^{***})$.
   \qed
  \section{Divergence--free\ $W$--valued random variables} \Label{domo}
   In this section we study the structure of\ $\domo{p}$. In classical
   analysis, zero--mean divergence--free vector fields generate rotations.
   On the other hand, the tangent space of the special orthogonal group\
   $SO(n)$\ can be identified with the space of skew symmetric\
   $n\mbox{\scriptsize $\times$}n$ matrices. Building on this correspondence,
   a 
   $W$--valued random variable $u$ was associated in~\cite{Z} with each
   sufficiently smooth random skew symmetric bounded operator $\AAA(\omega)$\
   in\ $H$,\ by
   \begin{equation} \Label {Atou}
    v=\sum_{i=1}^\infty \delta(\AAA e_i)e_i
   \end{equation}
   assuming the series converges in $L^p(\mu;W)$.
   Here\ $\calE\!=\!\{e_i\}_{_{i=1}}^{^\infty}$\! is a given smooth ONB.
   Under further smoothness and moment assumptions these $W$--valued random
   variables (vector fields) were then shown in~\cite{Z} to induce invariant
   flows in $W$, suggesting that they are, in our language, divergence--free.
%
   Note that~(\ref{Atou}) can now be written as $v\!=\!\ddelta\AAA^{\!\!^T}$.
   Indeed, for any\ $l\!\in\!W^*$,%
   \[ \clfs{v}{l}=\sum_{i=1}^\infty \delta(\AAA e_i)(l,e_i)_{_H}
                 =\delta\left(\sum_{i=1}^\infty \AAA e_i(l,e_i)_{_H}\right)
                 =\delta\left(\AAA\,\sum_{i=1}^\infty (l,e_i)_{_H}\,e_i\right)
                 =\delta(\AAA l)\]
   which by\ Lemma~\ref{weakversion} implies that $v\!=\!\delta \AAA^{\!\!^T}$,
   since\ $l\!\in\!W^*$\ was arbitrary.

   We show here (Theorem~\ref{domochar}) that {\it every} divergence-free
   $W$--valued $v$ can be obtained in this way, for some suitable skew
   symmetric random bilinear form\ $\AAA$\ on\ $W^*$.
  \begin{lemma} \Label{ueqminus}
   Let\ $u\!\in\!\domo{1}$\ \ be such that\
   $(\grad{u})^{\!^T}\!\!\in\!\Domdel{1}$. \ \  Then\
   $u=-\delta\,(\grad{u})^{\!^T}$.
  \end{lemma}
  \pf 
      For any given $G\!\in\!\calS(W^*)$\ apply Lemma~\ref {secmom} to obtain
      \[ E\clfs{u}{G}\!+\!E\clfs{\delta(\grad u)^T}{G}=0\]
      from which the result follows, since\ $G$\ was arbitrary.
  \qed
   \begin{theorem}  \Label{domochar}
    Let\ $v\!\in\!L^p(\mu\marker{;\Wss}),\ p\!\in\![1,\infty)$.\
    Then\ $v\!\in\!\domo{p}$\ \ iff there exists an\
    $\AAA\!\in\!\ddom_p\ddelta$\ with\ $\AAA\!+\!\AAA^T=0$\ such that
    $v\!=\!\ddelta \AAA$.
   \end{theorem}
   \pf
    Assume first that there exists an $\AAA$ as stated, and let
    $\Phi=\vph(\delta e_1,\ldots,\delta e_n)\!\in\!\calS_n$.
    Then, using the notation\ $\partial_x \Phi=\clf{\grad{\Phi}}{x}$\ for\
    $\Phi\!\in\!\calS$\ and\ $x\!\in\!W$,\ (and similarly\ $\partial^2_{xy}\Phi)$,
    \begin{eqnarray} \nonumber
     E\clfs{v}{\grad{\Phi}}
       &=&E\clf{\ddelta \AAA}{\sum_{i=1}^n
             \partial_{e_i}\vph(\delta e_1,\ldots,\delta e_n)e_i}\nonumber\\
       &\stackrel{\raisebox{0.15cm}{\footnotesize (\ref{weak})}}{=}&
             \sum_{i=1}^nE\partial_{e_i}\vph(\delta e_1,\ldots,\delta e_n)
                                               \delta(\AAA^Te_i)\nonumber\\
       &=&\sum_{i=1}^nE\clf{\AAA^Te_i}{\grad{\partial_{e_i}\vph}}\nonumber\\
       &=&E\sum_{i,j=1}^n\partial^2_{e_ie_j}\!\Phi\clfs{\AAA^Te_i}{e_j}
               \ =\ 0             \Label{orthog}
    \end{eqnarray}
    because ${\rm tr} AB\!=\!0$\ for all symmetric (respectively skew symmetric)\
    $n{\tiny \times}n$\ matrices\ $A$\ and\ $B$.\\
    Since\ $\Phi\!\in\!\calS$\ was arbitrary, it follows by definition
    that $\delta v$ exists and is $0$.

    For the converse, denote $u\!=\!(1+\calL)^{-1}v$ which also
    belongs to $\domo{p}$ by Lemma~\ref{domolemma}a). It then
    follows from Lemma~\ref{ueqminus} that
    \[ v-\calL(1+\calL)^{-1}v=(1+\calL)^{-1}v=-\delta(\grad(1+\calL)^{-1}v)^T.\]
    Recalling that $\calL\!=\!\delta\grad$, we thus have\ \ \
    $v=\delta\left(\grad(1+\calL)^{-1}v-(\grad(1+\calL)^{-1}v)^T\right)\ .$
   \qed
    \begin{remark} \Label{HodgeKodaira} \rm \small
     Proposition~\ref{vfdecomp} and Theorem~\ref{domochar} combined provide
     a unique decomposition of any element of\ $\dom_p\delta$\ as a sum of
     the gradient of a scalar random variable and the divergence of a
     random antisymmetric bilinear form. This constitutes an extension of
     I.~Shigekawa's first order Hodge--Kodaira theorey in Wiener space
     (cf. \cite{S}).
    \end{remark}
   \section{Wiener space \marker{valued} vector fields and
       \marker{their induced flows}}
   \Label{flowsection}
    One of the main interests in $W$--valued random variables $v$ is the
    possible existence of their generated quasiinvariant flows. More
    generally, we shall consider time--dependent vector fields and
    throughout, $I=[a,b]$ will be an compact interval in $\mathbb{R}$.
    \begin{definition} \Label{flowdef}
     A measurable mapping\
     $v:I\mbox{\footnotesize $\times$}W\longrightarrow W$ will be said
     to generate a flow\ $T=\{T_{s,t},\ s,t\!\in\!I\}$\ if $T$\ is
     jointly measurable and
    \begin{itemize}
     \item[\rm i)] \mbox{}\vspace{-1.8cm}

        \begin{equation}  \Label{solveseqn}
          T_{s,t}(\omega)=\omega+\int_s^t v_r(T_{s,r}(\omega))\,dr,
                 \hspace{1.5cm}\forall s,t\!\in\!I,\ \ \ \ {\rm a.s.}
        \end{equation}
      \item[\rm ii)]\mbox{}\vspace{-1.5cm}

        \begin{equation}  \Label{flow}
         T_{r,t}(\omega)=T_{s,t}(T_{r,s}(\omega))
               \hspace{1cm}\forall r,s,t\!\in\!I\ \ \ \ {\rm a.s.}
        \end{equation}
      \item[\rm iii)] \mbox{}\vspace{-1.5cm}

        \begin{equation}  \Label{abscont}
         T_{s,t} \mu \ll \mu, \quad \text{where\ } \quad
        \Lambda_{s,t} = dT_{s,t}^* \mu /d\mu
         \ \ \ \ \ \ \forall s,t\!\in\!I.
        \end{equation}
    \end{itemize}
    $T$ will be said to be the {\bf unique} flow generated by\ $v$\ if
    whenever $S=\{S_{s,t},\,s,t\!\in\!I\}$\ satisfies (i)--(iii) as well,
    \ \ \  $S_{s,t}=T_{s,t}$\ a.s. for all $s,t\!\in\!I$.
    \end{definition}

    In~\cite{Cru} A. Cruzeiro proved the existence of a flow for $H$--valued
    (time independent) vector fields whose gradient and divergence have
    finite exponential moments, and G. Peters noted in~\cite{P} that the
    weaker operator norm of\ $\grad{v}$\ could be used in the finite
    exponential moment assumption (instead of the usual Hilbert-Schmidt norm).
    In~\cite[Section 5.3]{UZ99} the result of Cruzeiro was extended and its
    assumptions relaxed.

    Based on Cameron--Martin's theorem on constant shifts and its subsequent
    generalizations, $H$--valued vector
    fields are indeed natural candidates to generate quasiinvariant flows.
    However, it is obvious that quasiinvariant shifts exist which do not
    necessarily act along $H$.

    In this section we present two general results. The first
    \marker{, Theorem~\ref{flowexists}, } states that
    ``representable" (to be defined) $W$--valued time dependent vector field
    $\{v_t,\ t\!\in\!I\}$, whose gradients and divergences satisfy some standard
    exponential moment conditions, generates a quasiinvariant flow.
    \marker{As such it extends~\cite{HUZ,Z}. }The
    decomposition of Proposition~\ref{vfdecomp} represent\ $v$'s\ $H$--shift
    and rotation components respectively.

    The second result,\marker{\ Theorem~\ref{musthavediv}, }shows that these
    are essentially all the $W$--valued random variables which do so, thus
    providing a qualitative description of the Wiener tangent space.
    \begin{definition} \Label{representable}
     Let\ $\calE\!=\!\{e_i\}_{i=1}^\infty$ be a given smooth ONB\ of $H$.
     \vspace{-.5cm}

     \begin{itemize}
      \item[a)] Denote\ \  $\pi_n(\om)\!=\!\sum_{i=1}^n \clf{e_i}{\om}e_i$,
        \ \ \ \ \  $W_n\!=\!\pi_n(W)$\ \ \ \ \ and\ \ \ \ \
        $\mu_n\!=\!\pi_n^*\mu$\\
        (without writing the dependence on\ $\calE$\ explicitly).
       Moreover, denote\ \ $\FF_n\!=\!\sigma(\pi_n)$\ \ and\ \
       $E_n(\,\cdot\,)\! =\!E(\,\cdot\,|\,\FF_n)$.\ \ A (not necessarily
       scalar) random variable on $(W,\calF,\mu)$\ will be said to
       be {\bf cylindrical} if it is $\FF_n$--measurable for some\
       $n\!\in\!\mathbb{N}$.
      \item[b)] A time-dependent vector field\
       ${\mathbf v}\!=\!\{v_t\}_{_{t\in I}}\!\in\!L^1({\rm leb}\times\mu;W)$\
       is said to be {\bf ${\bf \cal E}$--representable} if there exist\ \ \
       ${\mathbf v_i}\!=\!\{v_{i,t}\}_{_{t\in I}}\!\in\!L^1({\rm leb}\times\mu)
            ,\ i\!\in\!\mathbb{N}$,\ \ such that\ \
       $v_t\!=\!\sum_{i=1}^\infty v_{i,t}e_i$\ in\ $L^1({\rm leb}\times\mu;W)$.
     \end{itemize}
    \end{definition}
    \begin{theorem}  \Label{flowexists}
     Given a smooth ONB\ \  ${\cal E}\!=\!\{e_i\}_{i\in\mathbb{N}}$\ \ \ let
     ${\mathbf v}\!=\!\{v_t\}_{_{t\in I}}$\ be an\ ${\cal E}$-{\bf representable}
     vector field such that\  $v_t\!\in\!\dom_1\delta$\ \ for all $t\!\in\!I$,\
     and which moreover possesses a jointly measurable decomposition\
     $v_t\!\!=\!\!u_t\!+\!B_t,\ t\!\!\in\!I,$\ for which
     \vspace{-.5cm}

     \begin{itemize}
      \item[(i)] $u_t\!\in\!\mathbb{D}_{1,1}(H)$\ \ \ and\ \ \
                 $B_t\!\in\!\mathbb{D}^W_{1,1}\!(W)$\ with\
                 $\delta B_t\!=\!0$,\ \ \ for all
                 $t\!\in\!I$\\\hspace*{.2cm}
                 ($\mathbb{D}^{^W}_{1,1}\!(W)$\ is defined at the end of
                 subsection~\ref{Stochdiff})
      \item[(ii)] $\exists\theta\!>\!0$\ such that
      \vspace{-1.1cm}

       \begin{eqnarray} \Label{Gamma}
        \Gamma_{\!_H}(\theta)
          &:=&E\int_a^b\exp\theta\left\{\|\grad^{^{\!H}}u_t\|_{_{_{L(H)}}}
               \!+\!|\delta u_t|\right\}\,dt<\infty\hspace{1.5cm}{\rm and}
               \hspace{2cm}\\
        \Gamma_{\!_W} (\theta)     \Label{zakai-1}
          &:=&\sup_n E\int_a^b
             \exp\theta\|\pi_n\grad^{^W} B_t\|_{_{L(W)}}\;dt<\infty.
        \end{eqnarray}
      \end{itemize}
      Then ${\mathbf v}$ generates a flow\ $\{T_{s,t},\,s,t\!\in\!I\}$ with
      \begin{equation} \Label{flowrho}
       \frac{dT^*_{s,t}\mu}{d\mu}
                        =\exp\left\{\int_s^t\delta v_r(T_{t,r})\,dr\right\}
      \end{equation}
       and for all\ $p>1$\ and $|t\!-\!s|<\frac{\theta}{2p}$,
      \begin{equation} \Label{rhoflowbound}
       E\left(\frac{dT^*_{s,t}\mu}{d\mu}\right)^p
         \le e^{1/p}\left(1+\frac{2p-2}{\theta}\,
                           \sqrt{\Gamma_H(\theta)\Gamma_W(\theta)}\right)\ .
      \end{equation}
       If, in addition, the paths $t\to v_t$ are a.s. continuous
       then the flow $\{T_{s,t},\,s,t\!\in\!I\}$ is unique.
    \end{theorem}
 \begin{remarks}\mbox{}\\ \rm \small\vspace{-.7cm}

  \begin{itemize}
   \item[(a)] The decomposition assumed in Theorem~\ref{flowexists} is not
     unique. In particular, it is not necessarily the one provided by
     Proposition~\ref{vfdecomp}.
   \item[(b)] It follows from~\req{zakai-1} by Fatou's lemma that \marker{}
    \vspace{-.6cm}

     \begin{equation} \Label{zakai-2}
       E\int_a^b\exp\theta\|\grad^{^{\!W}} B_t\|_{_{L(W)}}\,dt
                  \le \Gamma_{W}(\theta)< \infty,
     \end{equation}
    however the proof below needs the less elegant assumption~\req{zakai-1}
    itself. If\ $\calE$ is a Schauder basis of $W$, then~\req{zakai-1}
    is also {\em implied} by~\req{zakai-2}\ (possibly with a smaller $\theta$),
    but in general this need not be the case.
   \item[(c)]~Equation~\req{Gamma} and \req{zakai-2} actually imply that
      $u_t\!\in\!\Dp{{}}(H)$\ \ and\ \ $B_t\!\in\!\Dp{W}(W)$\ \ for all
      $t\!\in\!I$\ and $p\!\ge\!1$. However the initial assumption in (i)
      for $p\!=\!1$\ was needed to give meaning to $\grad^{^{\!H}}u_t$\ and\
      $\grad^{^{\!W}}B_t$\ in the first place.
   \item[(d)]~It follows from~(\ref{flowrho}) that if\ $v_r$ is
      divergence free for all $r$\ then the flow is measure preserving.
  \end{itemize}
 \end{remarks}
  \vspace{.3cm}

 Before proceeding with the proof, let us recall some relevant results
 pertaining to jointly measurable time dependent vector field
 $\eta=\{\eta_t(x)\}_{_{t\in I}}$ on $W_n$ as stated, for example,
 in~\cite[Section~5.1]{UZ99}. Assume that\ $\eta$\ is locally
 integrable in the $t$ variable,\ $C^1$\ in the $x$ variable and that
 for some $\theta>0$
  \begin{equation}   \Label{Gammadef}
   \Gamma^{\eta}(\theta):= E_n\int_a^b\exp\left\{\theta
        \left(\left\|\grad \eta_t\right\|_{L(W_n)}
        \!+\!\left|\delta\eta_t\right|\right)\right\}\,dt<\infty.
    \end{equation}
  Then $\eta$ generates a quasiinvariant flow
  $\{X_{s,t},\,s,t\!\in\!I\}$ whose Radon-Nikodym derivative
    \[ R_{s,t}:=\frac{d\!\left(X_{s,t}^{^*}\mu_n\right)}{d\mu_n}
                 =\exp\int_s^t\left(\delta\eta_r\right)(X_{t,r})\,dr \]
    satisfies, for any $p\!>\!1$\ and for all\ $s,t\!\in\!I$\ for which\
    $|t\!-\!s|<\sfrac{\theta}{p}$,
    \begin{equation} \Label{rhoLpbound}
     E_n\left(R_{s,t}\right)^p
       \le e^{^{1/p}}\left(1+\sfrac{p-1}{\theta}\,\Gamma^\eta(\theta)\right)
    \end{equation}
    (This bound is essentially given in~\cite[Theorem~5.1.3]{UZ99}).
    Indeed, with $\frac{1}{p}\!+\!\frac{1}{q}\!=\!1$,
    \begin{eqnarray} \Label{Gronwall}
     E_n\left(R_{s,t}(x)\right)^p\!
       &=&E_n\left(R_{s,t}(X_{s,t}x)\right)^{p-1}
        =E_ne^{^{\mbox{\small $(p\!-\!1)\int_s^t\delta
                          \eta_r(X_{t,r}X_{s,t}x)\,dr$}}} \nonumber\\
       &\le&1+\frac{p}{\theta}\int_IE_ne^{^{\mbox{\small $\sfrac{\theta}{q}\,
                      |\delta\eta_r(x)|$}}}\,R_{s,r}(x)\,dr \nonumber \\
       &\le&\left(1+\frac{p-1}{\theta}\int_I
            E_ne^{^{\mbox{\small $\theta\,|\delta\eta_r(x)|$}}}dr\right)
             +\frac{1}{\theta}\int_s^tE_n\left(R_{s,r}(x)\right)^p\,dr\ .
    \end{eqnarray}
    The second row follows from\ \ \ \ $e^{^{\int_{_\alpha}^{^\beta}f}}\!\le\!
    1+\frac{1}{\gamma}\,\int_{_{\alpha\wedge\beta}}^{^{\alpha\vee\beta}}
    e^{^\gamma\!f}$\ \ \ \  whenever\ \ $|\beta\!-\!\alpha|\!<\!\gamma$\
    and\ $f(x)\!\ge\!0$\ \ \ \  (here \ \  $\gamma\!=\!\frac{\theta}{p}$),\ \
    while the third row follows from Young's inequality\ \
    $yz\!\le\!\frac{y^p}{p}\!+\!\frac{z^q}{q},\ \ y,z\!>\!0$.
    The estimate~(\ref{rhoLpbound}) now follows by by applying Gronwall's
    inequality to~(\ref{Gronwall}). Note that for the purposes
    of~(\ref{rhoLpbound}), $\Gamma^\eta(\theta)$ could have been defined
    in~(\ref{Gammadef}) without the gradient term.
  \bigskip

  \noindent
 {\bf Proof of Theorem~\ref{flowexists}}\ \ \
    Consider\ \ $v$'s\ \ finite dimensional ``projections"\ \ \
    $v^{(n)}_t\!=\!E_n(\pi_n(v_t))$,\ \ $t\!\in\!I$, which may be written
    as\ $v^{(n)}_t\!=\wtv^{(n)}_t\circ\pi_n$,\ for an appropriate
    $\wtv^{(n)}_t:W_n\to W_n$. The use of the notation\ $\grad$\ and\ $\delta$\
    below both in\ $(W,\mu)$\ and in\ $(W_n,\mu_n)$\ should cause no confusion.
    \begin{lemma} \Label{projgraddel}
     For all $t\!\in\!I$\ and\ $n\!\in\!\mathbb{N},\ \
       \wtv^{(n)}_t\!\in\!\mathbb{D}^{^{W_{_n}}}_{1,1}\!(W_n)$\ and
     \begin{eqnarray}
      \grad{\wtv^{(n)}_t}\circ\pi_n
         &=&\grad{v^{(n)}}=E_n\left(\pi_n(\nabla^{\!^H}u_t)\pi_n\right)
                      +E\left(\pi_n(\nabla^{\!^W}B_t)\pi_n\right)
                                                        ,\Label{projgrad}\\
      \delta{\wtv^{(n)}_t}\circ\pi_n&=&\delta v^{(n)}
                     =E_n\left(\delta v_t\right)\ \Label{projdel}
     \end{eqnarray}
    \end{lemma}
     {\bf Proof of the Lemma}:\ 
     Since~(\ref{projdel}) has already been obtained in
     Proposition~\ref{martingale}, we only need to  prove~(\ref{projgrad}).\
     Indeed, for\ $N>n$\ and\ $\varphi\!\in\!C^\infty_{\rm b}(\mathbb{R}^N)$\
     denote
     \begin{equation} \Label{smoothen}
      \wn{\varphi}(x_1,\ldots,x_n)=
        \int\!\!\cdots\!\!\int\hspace{-1.05cm}\raisebox{-3.1ex}
                       {\tiny $\mathbb{R}^{^{N\!-\!n}}$}\hspace{.3cm}
                \varphi(x_1,\ldots,x_{_N})d\mu_{_{N-n}}(x_{n+1},\ldots,x_{_N}).
     \end{equation}
     \vspace{-.8cm}

     \noindent
     If\ $\Phi=\varphi(\delta e_1,\ldots,\delta e_{_N})$\  then\ \
     $E_n\Phi=\wn{\varphi}(\delta e_1,\ldots,\delta e_n)\in\calS$,\ \ \
     and since clearly\ \ \ $\frac{\partial \wn{\varphi}}{\partial x_i}
     =\wn{\frac{\partial\varphi}{\partial x_i}}$\ \ for each\
     $i=1,\ldots,n$,\ \  it follows that
     \[ \grad E_n\Phi
        =\sum_{i=1}^n\frac{\partial \wn{\varphi}}{\partial x_i}
               (\delta e_1,\ldots,\delta e_n)e_i
        =\sum_{i=1}^N\wn{\frac{\partial \varphi}{\partial x_i}}
               \!(\delta e_1,\ldots,\delta e_n)e_i\pi_n
        = E_n(\grad\Phi)\pi_n.  \]
     For any linear subspace\ $Z$\ of\ $W$,\ separable Banach space\
     $Y$\ and $p\!\in\![1,\infty)$,\ this identity can be extended by
     linearity and density to
     \begin{equation} \Label{gradEn}
      \grad E_nF=E_n\left(\grad^{^{\!Z}}F\right)\pi_n\hspace{1.5cm}\forall F
                                                 \!\in\!\Dp{Z}(Y).
     \end{equation}
     On the other hand it is straightforward to check that
     \begin{equation}  \Label{gradpin}
      \grad^{^{\!Z}}\pi_nF=\pi_n\grad^{^{\!Z}}F
             \hspace{1.5cm}\forall F\!\in\!\Dp{Z}(Y)
     \end{equation}
     and~(\ref{projgrad}) is a direct consequence of~(\ref{gradEn})\
     and~(\ref{gradpin}).
     \qed
     \medskip

    \noindent
    Returning to the proof of the theorem, and recalling the
    definition~(\ref{Gammadef}),
 \begin{eqnarray} \Label{Gamman}
 \Gamma^{^{(\wtv^{(n)})}}(\sfrac{\theta}{2})
  &=&   E \int_a^b \exp \sfrac{\theta}{2}
      \Bigl(\| \nabla \wtv_t^{(n)} \|_{_{L(W)}} + |\delta
                 \wtv_t^{(n)}|\Bigr) \,dt \notag \allowdisplaybreaks\\
 &  \overset{\raisebox{1ex}{\scriptsize{Lemma~\ref{projgraddel}}}}{=}&
                                                        E \int_a^b \exp
 \sfrac{\theta}{2} \Bigl( \|E_n(\pi_n(\grad^{\!^W}v_t)\pi_n\|_{_{L(W)}}
                 + |E_n \delta v_t|\Bigr)\,dt\notag\allowdisplaybreaks\\
 &\le& E\int_a^b \exp \sfrac{\theta}{2} E_n \Bigl( \|\pi_n \grad^{\!^H}
         v_t \pi_n\|_{_{L(H)}} + \|\pi_n\grad^{\!^W} B_t \pi_n\|_{_{L(W)}}
                + |\delta v_t|\Bigr) \,dt\notag\allowdisplaybreaks\\
 &  \overset{\raisebox{1ex}{\scriptsize{Jensen}}}{\le}&
                                          E\int_a^b \exp \sfrac{\theta}{2}
 \Bigl( \|\grad^{\!^H} u_t\|_{_{L(H)}} + \| \pi_n \grad^{\!^W}
   B_t\pi_n\|_{_{L(W)}}+|\delta u_t|\Bigr)\,dt\notag\allowdisplaybreaks\\
     &\le& \left[E \int_a^b \exp \theta \Bigl(\|\nabla^{\!^H}u_t\|_{L(H)}
     +|\delta u_t|\Bigr) \,dt\right]^{\frac{1}{2}} \cdot \notag
                                                     \allowdisplaybreaks\\
 && \hspace{1cm}\cdot\left[\sup_n E \int_a^b \exp \theta \|\pi_n \grad^{\!^W}
      B_t \pi_n\| \,dt\right]^{\frac{1}{2}}\notag\allowdisplaybreaks\\
 &\le& \Bigl( \Gamma_{\!_H}(\theta) \Gamma_{\!_W}(\theta)\Bigr)^{\frac{1}{2}}\,.
 \end{eqnarray}

    \noindent
    In particular,\ \ $\forall t\!\in\!I,\ \
    \wtv^{(n)}_t\!\in\!\bigcap_{_{p\ge 1}}\!L^p(\mu_n;W_n)
    \,\subset\,C^\infty(W_n,W_n)$,\ \ \ by the Sobolev embedding theorem.
    Consequently,\ recalling the facts presented above in the finite
    dimensional setup,\ $\wtv^{(n)}_t$ generates a quasiinvariant flow
    \{$\wtTn_{s,t},\ s,t\!\in\!I\}$\  on $W_n$ satisfying
    \begin{equation}  \Label{fidiflow}
     \wtTn_{s,t}(x)
        =x+\int_s^t
        \wtv^{(n)}_r\left(\wtTn_{s,r}(x)\right)\,dr\ ,
         \hspace{1.5cm} x \in W_n
    \end{equation}
    whose Radon--Nikodym derivative
    \begin{equation} \Label{Rnn}
     \wtron_{s,t}:=\frac{d\left(\wtTn_{s,t}{}^{^*}\!\mu_n\right)}{d\mu_n}
                       =\exp\int_s^t\left(\delta\wtv^{(n)}_r\right)
                                 \circ\wtTn_{t,r}\,dr
    \end{equation}
    satisfies, for any $p\!>\!1$\ and for all\ $s,t\!\in\!I$\ for which\
    $|t\!-\!s|<\sfrac{\theta}{2p}$,
    \begin{equation} \Label{rhonLpbound}
     E_n\left(\wtron_{s,t}\right)^p\le e^{^{1/p}}
      \left(1+\sfrac{2(p-1)}{\theta}\,\sqrt{\Gamma_H(\theta)\Gamma_W(\theta)}\right).
    \end{equation}
    Here~(\ref{rhoLpbound}) was applied with $\eta\!=\!\wtron$\ and $\theta$\
    replaced by\ $\frac{\theta}{2}$,\ and making use of~(\ref{Gamman}).
    \medskip

    \noindent
    The flow\ $\wtTn_{s,t}$ can now be ``lifted" from \ $W_n$\ to\ $W$ by defining
    $$ T_{s,t}^{(n)} (\om) := \widetilde{T}_{s,t}^{(n)} (\pi_n w) +
         (\om-\pi_n \om). $$
    It is straightforward to verify that\ $T_{s,t}^{(n)}$\ satisfies
    the flow equation
 \begin{equation} \Label{Hncontrol}
   T_{s,t}^{(n)}(\om)=\om+\int_s^t v_r^{(n)}\Bigl(T_{s,r}^{(n)}
   (\om)\Bigr)\,dr.
 \end{equation}
 and that
 \begin{equation} \Label{liftedRN}
  \rho^{(n)}_{s,t}:=\frac{d\left({T_{s,t}^{(n)}}^*\!\mu\right)}{d\mu}
                   =\wtron_{s,t}\circ\pi_n
     =\exp\int_s^t\left(\delta v^{(n)}_r\right)\circ T^{(n)}_{t,r}\,dr
 \end{equation}
 The solution's construction will be completed by showing convergence of
 $T^{(n)}_{s,t}$ to a $W$--valued process $T_{s,t}$ which will be the required
 flow.
 \begin{proposition} \Label{approxprop}
  Let $\eta_t^{[0]},\eta_t^{[1]},\ \  t\in I$,\ \  be two $C^1$
  cylindrical time-dependent vector-fields on $W$ such that\
  $\Gamma_i(\theta):=\Gamma^{(\eta^{[i]})}(\theta)\!<\!\infty$,\  for some\
  $\theta\!>\!0$ \ and\ $i=0,1$\ \ \mbox{\rm \Large (}cf.~(\ref{Gammadef})$\
  \ \rule[0.1cm]{0.3cm} {0.01cm}\   \Gamma^{(\eta^{[i]})}(\theta)$\ should be
  understood to mean\ $\Gamma^{(\wteta^{[i]})}(\theta)$\ where\ $\eta^{[i]}\!
  =\!\wteta^{[i]}\!\circ\!\pi_n$\ for some $n\!=\!n(i)$\ \mbox{\rm \Large )}.
  \\ Let $\{L_{s,t}^{[i]}, s, t\in I\}$\ be the unique flow generated on $I$
  by $\eta^{[i]}$\ as above. Then for any $p>1$\ there exists a finite
  positive constant $c=c\left(p,\theta,\Gamma_0(\theta),\Gamma_1(\theta)\right)$,
  increasing in its third and fourth arguments, such that for any $s\!\in\!I$
  \begin{equation} \Label{pertu}
   E\sup_{\underset{|t-s|\le\frac{p-1}{2p}\theta}{t\in I}}
                     \left\|L_{s,t}^{[1]}-L_{s,t}^{[0]}\right\|_W
   \le c \left(\int_I E\left\|\eta_r^{[1]}-\eta_r^{[0]}\right\|_{_W}^p
                                                   \,dr\right)^{\frac{1}{p}}
  \end{equation}
 \end{proposition}
  {\bf Proof of the Proposition}:\
   For a fixed $p\!>\!1$ denote\
   $\Ds\!=\!\{t\!\in\!I,\ \ |t-s|\le\frac{p\!-\!1}{2p}\,\theta\}$.\
   Next, let\ \ $D_t\!=\!\eta_t^{[1]}\!-\!\eta_t^{[0]}$\ \ and for every\ \
   $\la\!\in\![0,1]$\ consider the interpolated vector field\ \ \
   $\eta_t\laa=\la\eta_t^{[1]}\!+\!(1\!-\!\la)\eta_t^{[0]}
   =\eta_t^{[0]}\!+\!\la D_t$.\ Note that by convexity
   \begin{equation}  \Label{Gammalambda}
     \Gamma^{^{\eta\laa}}\!(\theta)\le\Gamma_0(\theta)+\Gamma_1(\theta)
   \end{equation}
   for every $\la\!\in\![0,1]$. Now let\ $L\laa_{s,t}$\ be the flow $\eta\laa$
   generates on\ $I$, with induced Radon-Nikodym derivative $\rho_{s,t}\laa$,\
   so that setting $Z_{s,t}^{[\la]} = d \eta_{st}^{[\la]}/d\la$ yields
   \begin{equation} \Label{Ldiff}
      L_{s,t}^{[1]}-L_{s,t}^{[0]}=\int_0^1 Z_{s,t}\laa\,d\la
                                                \hspace{1.5cm}s,t\!\in\!I
   \end{equation}
   and it holds that
   \[ Z_{s,t}\laa(\om)=\int_s^t D_r\left(L_{s,t}\laa\om\right)\,dr +
         \int_s^t\grad\eta_r\laa(L_{s,r}\laa\om)\ \ Z_{s,r}\laa(\om)\ dr\ ,\]
   and thus by Gronwall's lemma, for all\ $t\!\in\!\Ds$
   \[\|Z_{s,t}\laa(\om)\|_{\!_W}\le\int_{\Ds}\|D_r(L_{s,r}\laa\om)\|_{\!_W}
       \,dr\ \ e^{^{\mbox{\small $\int_{\Ds}\|\grad\eta_r\laa(L_{s,r}\laa)\|
                                                   _{_{L(W)}}\,dr$}}}\ \ . \]

   \noindent
   Inserting this estimate in~(\ref{Ldiff}), and denoting\ \
   $p_0\!=\!\frac{p\!+\!1}{2},\ \ q_0\!=\!\frac{p\!+\!1}{p\!-\!1}\ \ \
   (\frac{1}{p_0}\!+\!\frac{1}{q_0}\!=\!1$),\ \ \ we obtain
   \begin{eqnarray*}
    \lefteqn{E\sup_{t\in\Delta_s}
           \|L_{s,t}^{[1]}(\om)-L_{s,t}^{[0]}(\om)\|_{\!_W}
      \le E\sup_{t\in\Delta_s}\int_0^1\|Z_{s,t}\laa(\om)\|_{\!_W}\,d\la}
\allowdisplaybreaks\\
      &\le& E
       \int_0^1 \left(\int_{\Ds}\|D_r(L_{s,r}\laa\om)\|_{\!_W}\,dr\
        e^{^{\mbox{\small $\int_{\Ds}\|\grad\eta_r\laa(L_{s,r}\laa\om)\|
                                _{_{L(W)}}\,dr$}}}\right)\,d\la
\allowdisplaybreaks\\
      &\overset{\raisebox{1ex}{\scriptsize{Jensen}}}{\le}&
           \int_0^1\left(E\int_{\Ds}\|D_r(L_{s,r}\laa\om)\|_{\!_W}\,dr
        \ \ \frac{1}{|\Ds|}\,\int_{\Ds} e^{\mbox{\small $\frac{p-1}{2p}\theta\,
          \|\grad\eta_r\laa(L_{s,r}\laa\om)\|_{_{L(W)}}$}}\,dr\right)\,d\la
 \allowdisplaybreaks\\
      &\le&\int_0^1E\left(\int_{\Ds}\|D_r(L_{s,r}\laa\om)\|^{^{p_0}}_{\!_W}
                                          \,dr\right)^{\frac{1}{p_0}}\,
        E\left(\int_{\Ds}e^{\mbox{\small $q_0\frac{p-1}{2p}\theta\,
               \|\grad\eta_r\laa(L_{s,r}\laa\om)\|_{_{L(W)}}$}}\,dr\right)
                                               ^{\frac{1}{q_0}}\,d\la\\
      &\le&\int_0^1\left(E\int_{\Ds}\|D_r(\om)\|^{^{p_0}}_{\!_W}\,
                          \rho_{s,r}\laa(\om)\,dr\right)^{\frac{1}{p_0}}\,
         \left(E\int_{\Ds}e^{\mbox{\small $\frac{p+1}{2p}\theta\,
         \|\grad\eta_r\laa(\om)\|
        _{_{L(W)}}$}}\,\rho_{s,r}\laa(\om)\,dr\right)^{\frac{1}{q_0}}\,d\la
\allowdisplaybreaks\\ \\
      &&\hspace{.5cm}\mbox{and apply\ H\"older's\ inequality twice with the
      conjugate\ pair}\ \sfrac{p}{p_0},\ \sfrac{p}{p\!-\!p_0}\ :
\allowdisplaybreaks\\
      &\le&\int_0^1\left(E\int_{\Ds}\|D_r(\om)\|^{^p}_{\!_W}\,
                          dr\right)^{\!\frac{1}{p}}\,
        \left(E\int_{\Ds}{\rho_{s,r}\laa(\om)}^{^{\frac{p}{p-p_0}}}\,dr\right)
                                             ^{\!\frac{p-p_0}{pp_0}}\ \cdot
\allowdisplaybreaks\\
                   &&\hspace{.8cm}\ \cdot\
          \left(E\int_{\Ds}e^{\mbox{\small $\theta\,
     \|\grad\eta_r\laa(\om)\|_{_{L(W)}}$}}\,dr\right)^{^{\!\frac{p_0}{q_0p}}}
         \!\left(E\int_{\Ds}{\rho_{s,r}\laa(\om)}^{^{\frac{p}{p-p_0}}}\,dr
                                               \right)^{\!\frac{p-p_0}{pq_0}}
          d\la\ .
   \end{eqnarray*}

   \noindent
   The product of the second and fourth factors in the integrand
   may be estimated using~(\ref{rhoLpbound})
   \begin{eqnarray*}
    \left(\int_{\Ds}E{\rho_{s,r}\laa}^{^{\frac{2p}{p-1}}}\,dr\right)^\frac{p-1}{2p}
      \!\!&\le&\!\!\left(|\Ds|e^{^{\frac{p-1}{2p}}}\left(1+\sfrac{p+1}{p-1}\,
         \sfrac{1}{\theta}\Gamma^{^{(\eta_r\laa)}}(\theta)\right)\right)
                                                           ^{\frac{p-1}{2p}}
                                                                                \\
      &\le&\left(\sfrac{1}{p}\,e^{^{\frac{p-1}{2p}}}\,
          \left(\theta(p\!-\!1)+(p\!+\!1)\Gamma^{^{(\eta_r\laa)}}(\theta)\right)
                                                     \right)^{\frac{p-1}{2p}}
   \end{eqnarray*}
   while the third factor is bounded above by\
   $\left(\Gamma^{^{(\eta_r\laa)}}(\theta)\right)^{\frac{p-1}{2p}}$\!\!.
   This, in conjunction with~(\ref{Gammalambda}), completes the
   proof of the Proposition.
   \qed\\
   \bigskip

   We now apply Proposition~\ref{approxprop} to the cylindrical vector
   fields  $v_t^{(n)}$ and the flows\ $T^{(n)}_{s,t}$\ they generate as
   in~(\ref{Hncontrol}). Fix an arbitrary $p\!>\!1$.\ By
   Corollary~\ref{martingalecor}\ \ $v^{(n)}$ is a Cauchy sequence in\
   $L^p(I\!\times\!W,{\rm leb}\!\times\!\mu\,;\,W)$; it thus follows from
   the proposition, applied with\ $\frac{\theta}{2}$,\  that there exists a\
   $\gamma\!=\!\gamma(p,\theta)\!>\!0$\ and a\ $W$--valued process\
   $\{T_{s,t};\ s,t\!\in\!I\mbox{\ and\ }|s-t|\!\le\!\gamma\}$\ \ such that
   \ for all $s\!\in\!I,\ \
   \lim_{_{n\to\infty}}\! \sup_{_{t\in I:|t-s|\le\gamma}}\!
   \|T_{s,t}\!-\!T^{(n)}_{s,t}\|_{_{_W}}\!\!=\!0$\ \ in  probability.
   (Note that the constant $c$ in~(\ref{pertu}) depends (monotonically) on
   the $\Gamma^{^{[\wtvn]}}\!\left(\frac{\theta}{2}\right)$'s, but these are
   uniformly bounded - cf.~(\ref{Gamman})\,). Furthermore, since we have
   almost sure uniform convergence along a subsequence, $T_{s,t}$ is almost
   surely continuous in\ $t$.

   \noindent
   For all $r,s,t\!\in\!I$\ any two of which are at a distance not larger
   than\ $\gamma$, the flow property~(\ref{flow}) obviously holds
   for $T^{(n)}$ and thus for $T$ as well.  This allows us
   to extend $T_{s,t}$ to all\ $s,t\!\in\!I$:\ \ \
   $T_{s,t}(\om)\!=\!T_{s_{m-1},s_m}(T_{s_{m-2},s_{m-1}}
                        (\cdots T_{s_1,s_2}(T_{s_0,s_1}(\om))\cdots))$
   for any sequence\ $s\!=\!s_0,s_1,\ldots,s_{m\!-\!1},s_m\!=\!t$\ in\ $I$\
   such that\ $|s_i\!-\!s_{i-1}|\!\le\!\gamma,\ \ 1\!\le\!i\!\le\!m$.
   This extension\ $\{T_{s,t},\ s,t\!\in\!I\}$\ is of course independent of
   the connecting sequence\ $\{s_i\}$,\ is a.s. continuous and satisfies the
   flow property~(\ref{flow}) as well.

   \noindent
   To show that\ $T_{s,t}$\ is quasiinvariant, fix any $p>1$\ and assume
   first that $|t\!-\!s|\!<\!\frac{p\!-\!1}{2p}\,\theta$.\ Since\
   $T_{s,t}^{{(n)}^*}\!\mu$\ and\ $\mu$\ are equivalent for every\
   $n\!\in\!\mathbb{N}$\ we need to verify that the sequence of respective
   Radon-Nikodym derivatives\ $\{\rho^{(n)}_{s,t}\}$ is uniformly integrable.
   However, this follows immediately from~(\ref{rhonLpbound}), and moreover
   by~\req{Rnn}
   \[\rho_{s,t}:=
       \frac{dT^{^*}_{\!s,t}\,\mu}{d\mu}=\lim_{n\to\infty}\rho\upn_{s,t}
               =\exp\int_s^t\left(\delta v_r\right)\circ T_{t,r}\,dr\]
   which is~\req{flowrho} (the last equality follows from~\req{liftedRN}
   and Lemma~\ref{martingale}). In addition, taking the limit as\ $n\to\infty$
   in~\req{rhonLpbound} and applying Fatou's lemma, the bound~\req{rhoflowbound}
   is obtained.
   The resulting quasiinvariance then holds for an arbitrary pair\
   $s,t\!\in\!I$ by connecting $s$ and $t$, if necessary, by a finite
   sequence\ $\{s_i\}$\ as above.
   \medskip

   We now proceed to show that\ $T_{s,t}$\ satisfies~(\ref{solveseqn}). This
   will be achieved as the limit in an appropriate sense of~\req{Hncontrol} as
   $n\!\to\!\infty$, the only nontrivial convergence being that of the integral
   term. Indeed, again fixing some\ $p\!>\!1$\ and assuming
   $|t\!-\!s|\!<\!\frac{p-1}{2p}\,\theta$,
   \begin{eqnarray} \Label{eqapprox}
    \int_s^t E\left\|v_r (T_{s,r}(\om))
              - v\upn_r(T\upn_{s,r}(\om))\right\|_{_W}\!\!dr
    &\le&\!\int_s^tE\left\|v_r(T_{s,r}(\om))
              \!-\!v_r(T\upn_{s,r}(\om))\nonumber\right\|_{_W}\!\!dr\hspace{3cm}\\
     &&\hspace*{-.35cm}+\int_s^tE\left\|v_r\left(T\upn_{s,r}(\om)\right)\!
                -\!v\upn_r\left(T\upn_{s,r}(\om)\right)\right\|_{_W}\!\!dr
   \end{eqnarray}
   \medskip
    Let\ $\eps\!>\!0$\ and $\frac{1}{p}\!+\!\frac{1}{q}\!=\!1$.
    Since
    \begin{eqnarray*}
     \lefteqn{E\int_s^t\left(\|v_r(T_{s,r}(\om))\|_{_W}+
                 \|v_r(T\upn_{s,r}(\om))\|_{_W}\right)\,dr}\hspace{2cm}\\
      &&=\int_s^tE\|v_r\|_{_W}\left(\rho_{s,r}+\rho\upn_{s,r}\right)\,dr\\
      &&\le\left(\left(\int_s^t\rho^q_{s,r}dr\right)^{^{\frac{1}{q}}}
         +\left(\int_s^t\rho^{{(n)}^q}_{s,r}dr\right)^{^{\frac{1}{q}}}\right)
                  \,\left(\int_s^tE\|v_r\|_{_W}^p\,dr\right)^{\frac{1}{p}}\\
      &&\le M\,\|{\mathbf v}\|_{_{L^p({\rm leb}\times\mu)}}<\infty
    \end{eqnarray*}
    where $M$ doesn't depend on $n$, there exists an $\eta\!>\!0$\ such that
    for every measurable\ $B\!\subset\![s,t]\times W$\ with\
    $({\rm leb}\!\times\!\mu)(B)\!\le\!\eta$
    \[ \underset{B}{\int\int}\left(\|v_r(T_{s,r}(\om))\|_{_W}+
       \|v_r(T\upn_{s,r}(\om))\|_{_W}\right)\,dr\,d\mu\le\eps.\]
    In particular choose the set $B$ provided by Lusin's theorem
    (\ \ $({\rm leb}\!\times\!\mu)(B)\!\le\!\eta$\ \ ) and the bounded
    continuous function $g$ on $[s,t]\!\times\!W$ for which
    $v{\mathbf 1}_{_{B^C}}=g{\mathbf 1}_{_{B^C}}$ a.e. Then, splitting
    the expectation over $B$ and $B^{\mbox{\tiny $C$}}$,
    \begin{eqnarray*}
     \lefteqn{\limsup_{n\to\infty}\int_s^tE\left\|v_r(T_{s,r}(\om))
        \!-\!v_r(T\upn_{s,r}(\om))\right\|_{_W}\!\!dr}\hspace{2cm}\\
         &&\le\ \ \ \eps\ +\lim_{n\to\infty}\int_s^tE\left\|g(T_{s,r}(\om))
          \!-\!g(T\upn_{s,r}(\om))\right\|_{_W}\!\!dr\ =\ \eps.
    \end{eqnarray*}
    Since\ $\eps$ was arbitrary, the first term in the right hand side
    of~\req{eqapprox} goes to zero. On the other hand, and referring
    to~\req{rhonLpbound}, the second term is bounded by
    \[ \int_s^t E\|v_r\!-\!v_r\upn\|_{_W}\rho\upn_{s,r}\,dr \le C\,
      \left(\int_s^t E\|v_r\!-\!v\upn_r\|^{^p}_{_W}\,dr\right)^{\frac{1}{p}}
                                  \underset{n\to\infty}{\longrightarrow}0\]
    with $C=\left((t\!-\!s)e^{^{1/q}}\left(1+\sfrac{2(q-1)}{\theta}\,
     \sqrt{\Gamma_H(\theta)\Gamma_W(\theta)}\right)\right)^{\frac{1}{q}}$.
    Thus for all $|t\!-\!s|\!<\!\frac{p-1}{2p}\,\theta$\  all the terms
    of~\req{Hncontrol} converge to those of \req{solveseqn}. The latter
    equation then holds as well for any arbitrary pair $s,t\!\in\!I$ as a
    result of $T$'s flow property which has already been proved.

    Finally to show uniqueness we first observe that since the vector field
    $v_t$ was assumed to possess continuous paths almost surely,
    any flow $S_{s,t}$ generated by it is a.s. continuously differentiable both
    in $t$ and in $s$. For such a flow
%
    define,\ \ for any fixed\ $s\!\in\!I$,
    \[    U_{s,t}=T_{t,s}\circ S_{s,t},\hspace{1.5cm}t\!\in\!I\, ,\]
    where $T$ is the particular flow constructed above. Our aim is to show
    that\ $\dot{U}_{s,t}=0$\ for all $t\!\in\!I$,\ \ a.s.
    (\,here $\dot{A}_t:=\frac{d}{dt}A_t$\,). Indeed
    \begin{equation} \Label{Udot}
      \dot{U}_{s,t}
           =\dot{T}_{t,s}(S_{t,s})+\grad^W\!T_{s,t}(S_{s,t})\,v_t(S_{s,t})
           =\left(\dot{T}_{t,s}+\grad^W\!T_{s,t}\,v_t\right)\circ S_{s,t}.
    \end{equation}
    For every fixed $t$, and by quasiinvariace, this expression will be $0$
    for one flow $S_{s,t}$ if it is so for any other. Since it {\em is} zero
    when\ $S_{s,t}\!=\!T_{s,t}$\ (in which case $U_{s,t}(\omega)\!=\!\omega$),
    we conclude that $\dot{U}_{s,t}=0$\ a.s., $\forall t\!\in\!I$,\ for {\it any}
    flow $S_{s,t}$. The quantifiers a.s. and $\forall t$ can now be interchanged
    since the processes in~(\ref{Udot}) have a.s. continuous paths.
    This completes the proof of Theorem~\ref{flowexists}.
%
%
    \qed
%
%
%

%
   \bigskip

   The next result shows that the existence of ${\mathbf v}$'s divergence
   in Theorem~\ref{flowexists} is ``nearly" a necessary assumption.
   \begin{theorem}   \Label{musthavediv}
    Let\ \  $\{v_t\}_{t\in I}:W\to W$\ \  be a time dependent vector field
    which generates a quasiinvariant flow\ $T_{s,t}$ with Radon-Nikodym
    derivative $\Lambda_{s,t}$.\
    Then\ \ $v_s\!\in\!\dom_1\delta$\ \ \ for any\ $s\!\in\!I$\ for which
    \begin{equation}  \Label{Lambdadot}
    \Lambda'_{s,s}:=\lim_{t\to s} \frac{\Lambda_{s,t}-1}{t-s}
    \end{equation}
     exists weakly in\ $L^1(\mu)$,\ in which case\ $\delta v_s=\Lambda'_{s,s}$.
     In particular, if $T_{s,t}$ is measure preserving on $I$ then
     $\delta v_s$\ exists and is zero for every\ $s\!\in\!I$.
   \end{theorem}
   \pf
    Let\ $\Phi\!\in\!\calS$\ be arbitrary. It is a direct consequence
    of ~(\ref{solveseqn}) that for every\ $s,t\!\in\!(a,b)$
    \[\Phi(X_{s,t})=\Phi(X_{s,s})+\int_s^t\clf{\grad{\Phi}(X_{s,r})}{v_r(X_{s,r})}\,dr
                                                     \hspace{1cm}{\rm a.s.}\]
    and thus
    \begin{eqnarray*}
      E\Lambda_{s,t}\Phi=E\Phi(X_{s,t})
         &=&E\Phi+\int_s^tE\clf{(\grad{\Phi})\circ X_{s,r}}{v_r\circ X_{s,r}}\,dr\\
         &=&E\Phi+\int_s^tE\Lambda_{s,r}\clf{\grad{\Phi}}{v_r}\,dr.
    \end{eqnarray*}
    It follows that\ $E\Lambda_{s,t}\Phi$ is differentiable in\ $t$\ and that\ \
    \[ E\Lambda'_{s,s}\Phi=\frac{\partial\{E\Lambda_{s,t}\Phi\}}{\partial t}
         \raisebox{-.2cm}{\rule{.1mm}{6.5mm}}_{\,t=s}=
E\clf{\grad{\Phi}}{v_s}
= E(\Phi \delta v_s)\,,
\]
    where the first equality results from the existence of~(\ref{Lambdadot}),
    thus proving the theorem.
   \qed

\section{Additional Remarks} \Label{addrem}
 In this final section we address two additional aspects of the flows
 introduced in Section~\ref{flowsection}: their adaptedness and their
 relevance to an associated PDE.
\subsection*{I.\ Adapted Flows:}
Let $(W, H, \mu)$ be an AWS and $\{\Pi^\theta, 0 \le \theta \le 1\}$ a
continuous, strictly increasing resolution of the identity on $H$.  Let
(\cite{UZ97}, 2.6 of~\cite{UZ99})
$$
\calF_\theta = \sigma \Bigl\{\delta (\Pi^\theta h), h \in H \Bigr\}, \quad 0 \le \theta
\le 1
$$
be the filtration induced by $\Pi^\bfcdot$ on
$(W, H, \mu)$.
In what follows we assume that $\Pi^\theta W^* \subset W^*$ for all
$\theta \in [0,1]$.  This can be easily verified for the classical Wiener
space.

\begin{definition}[\cite{UZ97} or Section~2.6 of~\cite{UZ99}]
\label{def-6.1} An $H$-valued random variable $u$ is
$\calF_\theta$ measurable if $(u,h)_H$ is $\calF_\theta$
measurable for all $h\in H$. Moreover, $u\in H$ is \emph{adapted}
if $(u, \Pi^\theta h)_H$ is $\calF_\theta$ measurable for all
$\theta \in [0,1]$.
\end{definition}

 \begin{definition} \Label{def-6.2}
 A\ $W$-valued random variable $u$ is adapted if there exists a sequence of
 $H$-valued r.v.'s $u_n$ $n=1,2,\dots$ which are adapted and
 $|u-u_n|_W\overset{\mathrm{a.s.}}{\longrightarrow}\:0$. Let $u$ be
 a $W$ or $W^{**}$ valued r.v.\ then $\Pi^\theta u \in W^{**}$ is
 defined by $$ \clfsss{\Pi^\theta u}{h} = \clf{\Pi^\theta h}{u}
 \hspace{1.5cm}
   \forall h\!\in\!W^*\,.$$
\end{definition}
 We prepare the following lemma for later reference.
 \begin{lemma} \Label{lem-6.1} (a)~If $U$ is a quasiinvariant adapted
  $W$--valued random variable and $\alpha$ is $\calF_\theta$--measurable,
  then $\alpha\circ U$ is also $\calF_\theta$--measurable.\\ \vspace{-.65cm}

 \noindent
 (b)~Let $U_1, U_2$ be as in (a) such that\ \
  $\Pi^\theta U_1\!=\!\Pi^\theta U_2$\ \ for some $\theta\!\in\![0,1]$.\
  Then if $\alpha$ is a $\calF_\theta$--measurable random variable,
  $a\circ U_1\!=\!\alpha\circ U_2$. Similarly, if $v$ is a $W$-valued adapted
  random variable, then $\Pi^\theta(v\circ U_1)\!=\!\Pi^\theta(v\circ U_1) $.
 \end{lemma}
 \pf Since $\alpha$ is $\calF_\theta$--measurable it is the a.s.\
 limit of polynomials in $\delta(\Pi^\theta h)$,\ $h\!\in\!W^*$\
 (recall that we have assumed $\Pi^\theta h\!\in\!W^*$), so that
 $\alpha\circ U$ is the a.s.\ limit of polynomials in\
 $\delta(\Pi^\theta h)\circ U = \clf{\Pi^\theta h}{U}$.
 Since $U$ is adapted, each $\clf{\Pi^\theta h}{U}$\ is
 $\calF_\theta$--measurable and thus so is $\alpha\circ U$. This proves~(a),
 and~(b) follows directly from~(a). \qed
 \begin{proposition} \Label{prop-6.1}
  Assume that $v_t(\om)$ is adapted for every $t\in I$\ and satisfies the
  conditions of Theorem~\ref{flowexists} (including the a.s. continuity in $t$).
  Let $v_t^\theta(\om):=\Pi^\theta v_t(\om)$ and assume that $v_t^\theta(\om)$
  satisfies these conditions as well for all $\theta \in [0,1]$.
  Then the solution to \eqref{solveseqn}:
  \[ T_{s,t}(\om)=\om+\int_s^t v_r(T_{s,r})\,dr\ \]
  is also adapted.
 \end{proposition}
 \pf Let $\widetilde{T}_{s,t}^\theta(\om)$ solve the equation
 \begin{equation} \label{page3}
  \widetilde{T}_{s,t}^\theta(\om) = \om + \int_s^t v_r^\theta\!\circ\!
               \widetilde{T}_{s,r}^\theta(\om)\,dr
 \end{equation}
 hence
 $$ \Pi^\theta \widetilde{T}_{s,t}^\theta(\om)=\Pi^\theta\om
      + \int_s^t v_r^\theta\!\circ\!\widetilde{T}_{s,r}^\theta(\om)\,dr $$
 and $\Pi^\theta\widetilde{T}_{s,t}^\theta(\om)$ is\ $\calF_\theta$-measurable
 since it is a measurable function of $\{v_r^\theta, r\!\in\![s,t]\}$ and
 $\Pi^\theta \om$.\\

 \noindent
 On the other hand, by \eqref{solveseqn}
 $$ \Pi^\theta T_{s,t}(\om)=\Pi^\theta\om
       + \int_s^t \Pi^\theta\Bigl(v_r\!\circ\!T_{s,r}(\om)\Bigr)\, dr \,. $$
 Set $T_{s,t}^\theta(\om)=\Pi^\theta T_{s,t}(\om)+(I-\Pi^\theta)\om$, then
 \begin{align} \label{page4}
  T_{s,t}^\theta \om=\om + \int_s^t (\Pi^\theta v_r)\!\circ\!T_{s,r}(\om)\,dr
       =\om + \int_s^t (\Pi^\theta v_r)\!\circ\!T_{s,r}^\theta(\om)\,dr\,.
 \end{align}
 Comparing \eqref{page3} with \eqref{page4} yields, by uniqueness,
 $\widetilde T_{s,t}^\theta(\om)=T_{s,t}^\theta(\om)$.  Hence
 $\lip T_{s,t}^\theta(\om),\Pi^\theta h\rip
         =\lip T_{s,t}(\om),\Pi^\theta h \rip$ is $\calF_\theta$-measurable
 and since $\theta \in [0,1]$ was arbitrary, $T_{s,t}(\om)$ is adapted.
 \qed

 \begin{remark}  \rm\small
  The flow considered by Cipriano and Cruzeiro \cite{CC} is of the type
  considered in this proposition.  Let $(W, H, \mu)$ be the classical
  $d$-dimensional Wiener process
  $$ \om = \left\{ \begin{pmatrix}\om'(t)\\\vdots\\\om^d(t)\end{pmatrix}
     , \quad t \ge 0 \right\} \quad \text{and} \quad \Bigl(v(\om)\Bigr)_\bfcdot
    =\begin{pmatrix}(v'(\om))_\bfcdot\\ \vdots\\ (v^d(\om))_\bfcdot
    \end{pmatrix}_\bfcdot\hspace{.5cm}.$$
  In the case considered in \cite{CC}
  $$ \Bigl(v^i(\om)\Bigr)_\bfcdot = \sum_{j=1}^d \int_0^\bfcdot
          a_i^j(\om, s) d \om^j(s) + \int_0^\bfcdot b^i (\om, s)\, ds $$
  and for every $1\!\le\!i,j\!\le\!d,\ s\!\ge\!0$,\ the coefficients\
  $a_j^i(\om, s)$ and $b^i(\om, s)$ are $\calF_s$ measurable. The
   assumptions in \cite{CC} guarantee that $v$ satisfies those of
  Proposition~\ref{prop-6.1}.
 \end{remark}

 \subsection*{II.\ The equation $\frac{df(t,\om)}{dt} = \delta (\AAA(\om)
  \cdot \grad f(t,\om)),\; (\AAA + \AAA^T=0)$.}

 Let $T_t\om$, $t\in R$, $T_0\om=\om$, be a measure preserving transformation
 on $W$.  $T_t\om$ is said to be a \textit{stationary process\/} if for any $n$,
 any $t_1, \dotsc, t_n$, any smooth $\varphi_i(\om)$ and any $\tau$
 \begin{equation}
 \label{star}
 \Law \Bigl\{ \varphi_1 (T_{t_1} \om), \dotsc, \varphi_n (T_{t_n} \om)\Bigr\}
 = \Law \Bigl\{ \varphi_i (T_{t_1+\tau} \om), \dotsc, \varphi_n (T_{t_n+\tau}
 \om) \Bigr\}
 \end{equation}
 A flow which is also a stationary process will be called a
 \textit{stationary flow}.  Note that if $T_t \om, t\in R$ is a measure
 preserving flow then it is also a stationary flow.

 \begin{proposition}
 \label{prop}
 Let $\AAA(\om)$ be a measurable and skew symmetric transformation on $H$.
 Further assume that $\AAA(\om)$ transforms $\DD_{p,1} (H)$ into
 $\DD_{p,1} (H)$.  Let $B(\om) = \sum_1^\infty \delta (\AAA (\om) e_i) e_i$,
 where $e_i, i=1,2, \dots$ is a smooth ONB on $H$, converges in
 $L^1(\mu; W)$, and assume that $T_t\om$, $t\in R$ solves
 $$
 \frac{dT_t\om}{dt} = B(T_t \om), \qquad T_0\om = \om
 $$
 and $T_t\om$, $t\in R$, is a \textit{stationary process}.
 Let $f_0(\om)$ be a smooth functional on $W$, for which
 $f_0 (T_t \om) \in \DD_{p,2}$.

 Then $f(t,\om): = f_0(T_t \om)$ solves the equation: 
 \begin{equation}
 \Label{p30} \frac{df(t,\om)}{dt} = \delta (\AAA (\om) \grad f(t,\om)\marker{)},
 \end{equation}
 $f(0, \om) = f_0 (\om)$.
 \end{proposition}

 \pf
 For any smooth $\varphi(\om)$ we have by stationarity
 \begin{align*}
 E \frac{1}{\eps} \Bigl[
 \Bigl(f_0 (T_{t+\eps} \om) - f_0 (T_t \om)\Bigr) \varphi (\om)\Bigr]
 & = \frac{1}{\eps} E \Bigl( f_0 (T_t\om) \cdot \varphi (T_{-\eps} \om) -
 f_0 (T_t\om) \varphi (\om) \Bigr) \\
 & = E f_0 (T_t\om) \frac{1}{\eps} \Bigl( \varphi (T_{-\eps } \om) -
 \varphi(\om) \Bigr)
 \end{align*}
 and since $(d \varphi (T_t \om)/dt)_{t=0}
  = \delta ( \AAA(\om \grad \varphi (\om) ) $
 (cf. \cite{HUZ} or eqn.~1.10 of \cite{Z})
 $$
 E \varphi(\om) \frac{d f_0 (T_t \om)}{dt} = -
 E\Bigl( f_0 (T_t \om) \delta (\AAA (\om) \grad \varphi (\om))\Bigr)
 $$
 integrating by parts yields
 \begin{align*}
 E \varphi (\om) \frac{d f_0 (T_t \om)}{dt}
 & = - E \Bigl( \grad (f_0 (T_t \om)), \AAA \grad \varphi\Bigr)\\
 & = E\Bigl( \delta (\AAA \grad f_0 (T_t \om)) \varphi (\om) \Bigr)
 \end{align*}
 and \eqref{p30}  follows. \qed

 \begin{corollary}
 If in addition to the assumptions of proposition~\ref{prop}, $
 \grad f (t,\om) \in W^*$ a.s.\ for every $t$, then
 \begin{equation}
 \label{starstar}
 \frac{df(t,\om)}{dt} =
 \clf{\grad f(t,\om)}{B(\om)}  \qquad
 f(0,\om) = f_0 (\om)
 \end{equation}
 \end{corollary}

 \pf
 Let $\grad f(t,\om) = \sum_i \psi_i (t,\om) e_i$, then since
 $f_0 (t,\om) \in \DD_{p,2}$,
 $\psi_i (t,\om) \in \DD_{p,1}$ and
 \begin{align*}
 \delta (\AAA \grad f) & = \sum_1^\infty \delta \Bigl(\psi_i (t,\om) \AAA (\om)
 e_i\Bigr)
 \\
 & = \sum_1^\infty \psi_i (t,\om) \delta \Bigl(\AAA (\om) e_i\Bigr)
 - \sum_1^\infty \grad \psi_i (t,\om), \AAA (\om) e_i\\
 & = \lip B(\om), \grad f(t,\om)\rip - \sum_i \sum_j \grad_{e_i, e_j}^2
 f(t,\om) (e_j, \AAA e_i)\\
 & = \lip B(\om), \grad f(t,\om)\rip - 0\,.
 \end{align*}
 \qed

 In the converse direction
 \begin{proposition}
 Let $\AAA(\om)$ be a measurable and skew symmetric transformation on $H$
 transforming $\DD_{2,1}(H)$ into $\DD_{2,1}(H)$.  Assume that
 $f^j (t,\om) \in \DD_{p,2}$ and
 $\grad f^j \in W^*$, $j=1,2, \dots$ solves equation~\eqref{p30} for all
 $t$ and for the initial condition
  $f^j (0,\om) = \clf{e_j}{\om}$ where
 $e_j, j=1,2, \dots$ is a smooth ONB on $H$.  Then
 \begin{itemize}
 \item[(a)]
 $\sum_1^n \beta_j f^j (t+\tau_j, \om)$ solves \eqref{p30} under the
 initial condition $\sum_1^n \beta_j f^j (\tau_j, \om)$ for any
 $\beta_j, \tau_j$.

 \item[(b)]
 Let $\psi^n (t, \om) = \exp i \sum_1^n f^j (t+\tau_j, \om)$, then
 $\psi^n (t,\om)$ solves \eqref{p30} under the initial condition $$
 \psi^n(0, \om) = \exp i \sum_1^n f^j (\tau_j, \om)\,. $$
 \item[(c)]
  $T_t\om$ defined by
 $$
 T_t \om = \sum_{j=1}^\infty f^j (t, \om) e_j\,,
 $$
 is a stationary process.
 \item[(d)]
 If moreover,\eqref{p30} possesses a unique solution for every initial
 $f_0(\om) \in \DD_{p,2}$, then $T_s(T_t\om) = T_{s+t} \om$.
 \end{itemize}
 \end{proposition}

 \pf
 (a)~is trivial.
 (b)~since $f^j(t,\om) \in \DD_{p,2}$ and
 $\AAA + \AAA^T = 0$
 \begin{align}
 \label{secondstarstar}
 \delta\Bigl(\AAA \grad f^k (t, \om)\Bigr)
 & = \sum_j \delta ( \grad_{e_j} f^k (t,\om) \AAA e_j) \notag \\
 & = \sum_j \grad_{e_j} f^k (t,\om) \delta (\AAA e_j)
 - \sum_j\sum_k \grad_{e_i, e_j}^2 f^k (t,\om) (e_k, \AAA e_j)\notag\\
 &= \clf{\grad f^k (t,\om)}{\delta\AAA}\,.
 \end{align}
 On the other hand, differentiating $\psi^n(t,\om)$ with respect to $t$ and
 applying \eqref{secondstarstar} yields
 \begin{align*}
 \frac{d\psi^n (t,\om)}{dt} & = \psi^n (t,\om) \clf{\sum_{k=1}^n \grad
 f^k(t+\tau_j, \om)}{\delta \AAA}\\
 & = \clf{\grad \psi^n (t,\om)}{\delta\AAA} \\
 & = \delta \Bigl(\AAA \nabla \psi^n (t,\om)\Bigr)
 \end{align*}
 proving (b).  Turning to (c), by the last equation for $\tau_j=0$,
 $j=1,2,\dotsc,$
 \begin{equation}
 \label{starstarstar}
 E \psi^n (t, \om) = E \psi^n (0, \om)
 \end{equation}
 and since $e_i$ is an ONB, it follows that $f^i(t,\om)$ are
 $N(0,1)$ i.i.d.\ random variables; hence by the Ito-Nisio theorem
 $T_t\om$ is a measure invariant transformation.  Moreover, $T_t\om$ is a
 stationary process since \eqref{starstarstar} holds for any $t$ and any
 $\tau_1, \tau_2,
 \dots$.

 To show (d) note first that $f^i(t,\om) = \lip T_t\om, e_i\rip $ solves
 \eqref{p30} under $f^i(0,\om) = \lip \om, e_i\rip$, hence
 $\lip T_{t+\tau} \om, e_i\rip$ solves the same equation under
 $f_0^i(\om) = \lip T_\tau \om, e_i\rip$.
 On the other hand, as in (b), for any smooth
 $f_0(\om) = \widetilde{f}_0 (\lip\om, e_1\rip, \dots \lip\om, e_n\rip)$
 $f_0(T_t\om) = \widetilde{f}_0 (\lip T_t \om, e_1 \rip, \dots \lip T_t \om,
 e_n\rip )$ solves \eqref{p30}.
 In particular set $f_0^i(\om) = \lip T_\tau \om, e_i\rip$, then
 \begin{align*}
 f^i (t, \om) & = \lip T_\tau \om, e_i \rip\, \circ\, T_t \om\\
 & = \lip T_\tau (T_t \om) \rip, e_i \rip
 \end{align*}
 and (d) follows by the uniqueness of the  solution to \eqref{p30}.
 \qed

 

\begin{thebibliography}{99}
 \bibitem{BM} V. Bogachev and E. Mayer-Wolf, Absolutely continuous flows
    generated by Sobolev class vector fields in finite and infinite
    dimensions, {\sl J. Func. Anal. \bf 167} (1999), 1--68.
 \bibitem{CC} F. Cipriano and A.B. Cruzeiro, Flows associated to tangent
    processes on Wiener space, {\sl J. Funct. Anal. \bf 166} (1999), 310--331.
 \bibitem{Cru}A.B. Cruzeiro, \'{E}quations differentielles sur l'espace de
   Wiener et formules de Cameron--Martin non lin\'{e}aires, {\sl J. Funct.
   Anal. \bf 54} (1983), 206--227.
\bibitem{CM96} A.B. Cruzeiro and P. Malliavin, Renormalized differential
   geometry on path space, structural equation, curvature, {\sl J. Funct.
   Anal. \bf 139} (1996), 119--181.
 \bibitem{CM00} A.B. Cruzeiro and P. Malliavin, A class of anticipative
    tangent processes  on the Wiener space, {\sl C.R. Acad. Sci. Paris,
    \bf 333(1)} (2001), 353--358.
 \bibitem{D92} B.K. Driver, A Cameron-Martin type quasi-invariance theorem
    for Brownian motion on a compact manifold, {\sl J. Funct. Anal. \bf 110}
    (1992), 272--376.
 \bibitem{D95} B.K. Driver, Towards calculus and geometry on path spaces,
    {\sl Symp. Proc. Pure Math. \bf 57} (1995), 405--422.
\bibitem{FP} D. Feyel and de La Pradelle, Espaces de Sobolev gaussiens,
{\sl Ann. Inst. Fourier \bf 41}, (1991), 49--76.
 \bibitem{HUZ} Y. Hu, A.S. \"Ust\"unel and M. Zakai, Tangent processes on
    Wiener space, {\sl J. Funct. Anal. \bf 192} (2002), 234--270.
 \bibitem{K2} S. Kusuoka, Nonlinear transformations containing rotation and
    Gaussian measure. {\sl J. Math, Sci. Univ. Tokyo, \bf 10}, (2003), 1--40.
 \bibitem{L} P-L. Lions, Sur les \'equations diff\'erentielles ordinaires et
    les \'equations de transport, {\sl C.R. Acad. Sci Paris \bf 326(1)},
    (1998) 833--838.
 \bibitem{LT} M. Ledoux and M. Talagland, {\sl Probability in Banach spaces},
    Springer 1991.
 \bibitem{M} P. Malliavin, {\sl Stochastic Analysis}, Springer-Verlag,
    Berlin/New York, 1997.
\bibitem{MD} P. Malliavin and D. Nualart, Quasi sure analysis of stochastic
    flows and Banach space valued smooth functionals on the Wiener space,
    {\sl J. Funct. Anal. \bf 112} (1993), 287--317.
\bibitem{NZ} D. Nualart and M. Zakai, A summary of some identities of the
Malliavin calculus, In Stochastics Partial Differential Equations and
Applications II, G. Da Prato and L.~Tubero, editors.  {\sl Lecture Notes in
Mathematics \bf 1390}, 192--196, Springer 1989.
 \bibitem{P} G. Peters, Anticipating flows on the Wiener space generated by
    vector fields of low regularity, {\sl J. Funct. Anal. \bf 142} (1996)
    129--192.
 \bibitem{S} I. Shigekawa, De Rham-Hodge-Kodaira's decomposition on an abstract
    Wiener space, {\sl J. Math. Kyoto Univ. \bf 26}, 191--202 (1986).
\bibitem{S94} I. Shigekawa, Sobolev spaces of Banach-valued functions
associated with a Markov process, {\sl Prob. Th. Related Fields, \bf 99}
(1994) 425--441.
\bibitem{U} A.S. \"Ust\"unel,  An Introduction to Analysis of Wiener Space,
{\sl Lecture Notes in Math. \bf 1610}, Springer 1996.
 \bibitem{UZ97} A.S. \"Ust\"unel and M. Zakai, The construction of filtrations
    on abstract Wiener space, {\sl J. Funct. Anal. \bf 143} (1997) 10--32.
 \bibitem{UZ99} A.S. \"Ust\"unel and M. Zakai, {\sl Transformation of Measure
    on Wiener Space}, Springer-Verlag, New York/Berlin, 1999.
 \bibitem{Z} M. Zakai, Rotations and tangent processes on Wiener space,
    {\sl Seminaire de Probabilities}  {\small XXXVIII} 2004 to appear.
    (arXiv:math. PR/0301351).
 \end{thebibliography}
\end{document}